\documentclass[ptrf,preprint]{imsart}

\RequirePackage{amsthm,amsmath}
\RequirePackage[numbers]{natbib}
\RequirePackage[colorlinks,citecolor=red,urlcolor=blue,linkcolor=blue]{hyperref}
\RequirePackage{url,xspace,amssymb}
\RequirePackage{amssymb} \RequirePackage{amsfonts}
\RequirePackage{latexsym} \RequirePackage[utf8]{inputenc}
\RequirePackage{graphicx} \RequirePackage[english,french]{babel}
\RequirePackage{graphics}

\newcommand{\epr}{\hfill\hbox{\hskip 2pt
                \vrule width 5pt height 6pt depth 1.5pt}\vspace{0.1cm}\par}

\newcommand{\ve}{\varepsilon}

\newcommand{\cA}{{\cal A}}
\newcommand{\cB}{{\cal B}}
\newcommand{\cC}{{\cal C}}
\newcommand{\cD}{{\cal D}}
\newcommand{\cE}{{\cal E}}
\newcommand{\cF}{{\cal F}}
\newcommand{\cG}{{\cal G}}
\newcommand{\cH}{{\cal H}}
\newcommand{\cI}{{\cal I}}

\newcommand{\cK}{{\cal K}}

\newcommand{\cM}{{\cal M}}
\newcommand{\cN}{{\cal N}}

\newcommand{\cP}{{\cal P}}
\newcommand{\cQ}{{\cal Q}}

\newcommand{\cV}{{\cal V}}

\newcommand{\cY}{{\cal Y}}

\newcommand{\bA}{\mathbb A}

\newcommand{\bE}{\mathbb E}

\newcommand{\bH}{\mathbb H}
\newcommand{\bI}{{\mathbb I}}

\newcommand{\bN}{{\mathbb N}}

\newcommand{\bP}{{\mathbb P}}

\newcommand{\bR}{{\mathbb R}}

\newcommand{\mB}{\mathfrak{B}}

\startlocaldefs \numberwithin{equation}{section}
\theoremstyle{plain}
\newtheorem{theorem}{Theorem}[section]
\newtheorem{proposition}{Proposition}
\newtheorem{lemma}{Lemma}
\newtheorem{definition}{Definition}
\newtheorem{remark}{Remark}

\newtheorem{assertion}{Assertion}

\endlocaldefs

\begin{document}

\begin{frontmatter}
\title{Minimax and minimax adaptive estimation in multiplicative regression: locally bayesian approach}
\runtitle{Locally bayesian approach}

\begin{aug}
\author{\fnms{M.} \snm{Chichignoud} \ead[label=e1]{chichign@cmi.univ-mrs.fr}}

\runauthor{M. Chichignoud}

\affiliation{Université Aix-Marseille 1}

\address{Université Aix-Marseille 1\\
  LATP, 39 rue Joliot Curie, 13453 Marseille cedex 13, FRANCE,\\
  \printead{e1}\\
\phantom{E-mail:\ }}
\end{aug}
\selectlanguage{english}
\begin{abstract}
The paper deals with the non-parametric estimation in the regression
with the multiplicative noise. Using the local polynomial fitting
and the bayesian approach, we construct the minimax on isotropic
Hölder class estimator. Next applying Lepski's method, we propose
the estimator which is optimally adaptive over the collection of
isotropic Hölder classes. To prove the optimality of the proposed
procedure we establish, in particular, the exponential inequality
for the deviation of locally bayesian estimator from the parameter
to be estimated. These theoretical results are illustrated by
simulation study.
\end{abstract}



\end{frontmatter}

\section{Introduction}Let statistical experiment be generated by the couples of observations
$Y^{(n)}=(X_i,Y_i)_{i=1,...n},\:n\in\bN^{*}$ where $(X_i,Y_i)$
satisfies the equation
\begin{equation}
\label{model} Y_i=f(X_i)\times U_i,\quad i=1,\ldots, n.
\end{equation}
Here $f:[0,1]^d\rightarrow\mathbb{R}$ is unknown function and we are interested in estimating $f$ at a given point $y\in[0,1]^d$ from observation
$Y^{(n)}$.

The random variables (noise) $(U_i)_{i\in{1,\ldots,n}}$  are supposed to be independent and  uniformly distributed on $[0,1]$.

The design points $(X_i)_{i\in{1,...,n}}$ are deterministic and
without loss of generality we will assume that
$$
X_i\in \left\{1/n^{1/d},2/n^{1/d},\ldots, 1\right\}^{d},\:\:
i=1,\ldots,n.
$$
Along  the paper the unknown function $f$  is supposed to be smooth,
in particular, it belongs to the Hölder ball of functions
$\mathbb{H}_d({\beta},L,M)$ (see Definition \ref{def_holder} below).
Here ${\beta}>0$ is the smoothness of $f$, $M$ is the sum of upper bounds of
$f$ and its partial derivatives and $L>0$ is Lipschitz constant.

Moreover, we will consider only the functions $f$ separated away
from zero by some positive constant. Thus, from now on we will suppose
that there exists $0<A<M$ such that $f\in
\mathbb{H}_d({\beta},L,M,A)$, where
$$
 \mathbb{H}_d({\beta},L,M,A)=\left\{g\in\mathbb{H}_d({\beta},L,M): \inf_{x\in [0,1]^d}g(x)\geq A\right\}.
$$

\paragraph{Motivation.}
The theoretical interest to the multiplicative regression model (\ref{model}) with discontinuous noise
is dictated by the following fact.
The typical  approach to the study of the models with multiplicative
noise consists  in their transformation into the model with an additive noise
and in the application, after that, the  linear smoothing technique, based
on standard methods like kernel smoothing, local polynomials etc.
Let us illustrate the latter approach by the consideration of one of the
most popular non-parametric model namely multiplicative gaussian regression
\begin{equation}
\label{model-gauss} Y_i=\sigma(X_i) \xi_i,\quad i=1,\ldots, n.
\end{equation}
Here $\xi_i,\: i=1,\ldots, n $ are i.i.d. standard gaussian random variables and  the goal is to estimate
the variance $\sigma^{2}(\cdot)$.

Putting $Y^\prime_i=Y^2_i$ and $\eta_i=\xi_i^2-1$ one can transform the model (\ref{model-gauss}) into
the heteroscedastic additive regression:
$$
Y^\prime_i=\sigma^{2}(X_i)+ \sigma^{2}(X_i) \eta_i,\quad i=1,\ldots, n,
$$
where, obviously, $\bE\eta_i=0$. Applying any of the linear methods mentioned above to the estimation of $\sigma^{2}(\cdot)$ one can construct an
estimator whose estimation accuracy is given by $n^{-\frac{\beta}{2\beta+d}}$ and which is optimal in minimax sense (See Definition
\ref{def_minimax}). The latter result is proved under assumptions on $\sigma^{2}(\cdot)$ which are similar to the assumption imposed on the function
$f(\cdot)$. In particular, $\beta$ denotes the regularity of the function $\sigma^{2}(\cdot)$. The same result can be obtained  for any noise
variables $\xi_i$ with known, continuously  differentiable density, possessing
sufficiently many moments.
\par

The situation changes dramatically when one considers the noise with discontinuous distribution density. Although, the transformation of the original 
multiplicative model to the additive one is still possible, in particular, the model (\ref{model}) can be rewritten as
$$
Y^\prime_i=f(X_i)+ f(X_i) \eta_i,\quad Y^\prime_i=2Y_i,\:\:\eta_i=2u_i-1,\quad i=1,\ldots, n,
$$
the linear methods are not optimal anymore. As it is proved in Theorem
\ref{th:lower_bound} the optimal accuracy is given by $n^{-\frac{\beta}{\beta+d}}$. To achieve this rate the non-linear estimation procedure,
based on locally bayesian approach, is proposed in Section \ref{sectionMinimax}.

\par

Another interesting feature is the selection from given family of estimators (see \cite{Barron_Birge_Massart99}, \cite{Goldenshluger_Lepski08}). Such
selections are used for construction of data-driven (adaptive) procedures. In this context, several approaches to the selection from the family of
linear estimators were recently proposed, see for instance \cite{Goldenshluger_Lepski08}, \cite{Goldenshluger_Lepski09},
\cite{Juditsky_Lepski_Tsybakov09} and the
references therein. However, these methods are heavily based on the linearity property. As we already mentioned
the locally bayesian estimators are non-linear and in Section \ref{sectionAdaptive} we propose the selection rule from this family. It requires, in
particular, to develop new non-asymptotical exponential inequalities, which may have an independent interest.

\medskip

Besides the theoretical interest, the multiplicative regression model is applied in various domains, in particular,
in the image processing, for example, in so-called {\it nonparametric frontier
model} (see \cite{Allon_Beenstock_Hackman_Passy_Shapiro}, \cite{Simar_Wilson00}) can be considered as the particular case of the model (\ref{model}).
Indeed, the reconstruction  of the regression function $f$ can be viewed as the estimation of a production set $\cP$. Indeed,  $Y_i\leq
f(X_i),\:\forall i$, and, therefore, the estimation of $f$ is reduced to finding the upper boundary of $\cP$. In this context, one can also cite
\cite{Korostelev_Tsybakov93} dealing with the estimation of function's support. It is worth to mention that although nonparametric estimation in the
latter models is studied, the problem of adaptive estimation was not considered in the literature.


\paragraph{Minimax estimation.}The first part of the paper is devoted to the minimax over
$\bH_d(\beta,L,M,A)$ estimation. This means, in particular, that the
parameters $\beta,L,M$ and $A$ are supposed to be known {\it a
priori}. We find the {\it minimax rate of convergence}
(\ref{lower_sequence})  on $\bH_d({\beta},L,M,A)$ and propose the
estimator being optimal in minimax sense (see Definition
\ref{def_minimax}). Our first result (Theorem \ref{th:lower_bound}) in this direction consists in
establishing a lower bound for maximal risk on
$\bH_d({\beta},L,M,A)$. We show that for any
${\beta}\in\bR^{*}_{+},$ the minimax rate of convergence is bounded
from below by the sequence
 \begin{equation}
 \label{lower_sequence}
 \varphi_n({\beta})=n^{-\frac{{\beta}}{{\beta}+d}}.
 \end{equation}

Next, we propose the minimax estimator, i.e. the estimator attaining
the normalizing sequence (\ref{lower_sequence}). To construct the
minimax estimator we use so-called {\it locally bayesian estimation
construction}
which consists in the following. Let
$$
V_h(y)=\bigotimes_{j=1}^d\big[y_j-h/2,y_j+h/2\big],
$$
be the neighborhood around $y$ such that $V_h(y)\subseteq[0,1]^d$,
where $h\in (0,1)$ is a given scalar. Fix an integer number $b>0$ and let
$$
D_b=\sum_{m=0}^{b}\binom{m+d-1}{d-1}.
$$



Let $ \cP_b=\left\{p=(p_1,\dots,p_d)\::\:p_i\in\bN,\:0\leq |p|\leq b\right\},\:|p|=p_1+\dots+p_d $, we define the local polynomial
\begin{eqnarray}\label{polynome}
f_t(x)&=&\sum_{p\in\cP_b} t_p\left(\frac{x-y}{h}\right)^p\bI_{V_h(y)}(x), \quad
x\in\bR^d,\:t=(t_p \::\: p\in\cP_b),
 \end{eqnarray}
where $ z^p=z_1^{p_1}\cdots z_d^{p_d} $ for $ z=(z_1,\dots,z_d) $ and $\bI$ denotes the indicator function.  The local polynomial
$f_t$ can be viewed as an approximation of the regression function
$f$ inside of the neighborhood $V_h$ and $ D_b $ the number of coefficients of this polynomial. Introduce the following subset
of $\bR^{D_b}$
\begin{eqnarray}
\label{ensemble} \Theta\big(A,M\big)=\left\{t\in\bR^{D_b}:\:\:
2t_{0,...,0}-\|t\|_1\geq A,\:\: \|t\|_1\leq M\right\},
\end{eqnarray}
where $\|.\|_1$ is  $l_1$-norm on $\bR^{D_b}$. $ \Theta(A,M) $ can be viewed as the set of coefficients $ t $ such that $ A\leq f_t(x)\leq M $ for
all $ t\in\Theta(A,M) $ and for all $x$ in the neighbourhood $ V_h(y) $. Consider the {\it pseudo likelihood ratio}
$$
L_h\big(t,Y^{(n)}\big)=\prod_{i:\:X_i\in
V_h(y)}\big[f_t(X_i)\big]^{-1}\bI_{\big[0,
f_t(X_i)\big]}\big(Y_i\big),\quad t\in \Theta\big(A,M\big).
$$
Set also
\begin{equation}
 \label{critere}
    \pi_h(t)=\int_{\Theta(A,M)}\|t-u\|_1\: L_h\big(u,Y^{(n)}\big)du,\quad
    t\in \Theta\big(A,M\big).
 \end{equation}
  Let $\hat{\theta}(h)$ be the solution of the following minimization problem:
\begin{equation}
 \label{bayes-estimator}
\hat{\theta}(h)=\arg\min_{t\in \Theta(A,M)}\pi_h(t).
\end{equation}
The {\it locally bayesian estimator} $\bar{f}^h(y)$  of $f(y)$ is defined now as  $ \bar{f}^h(y)=\hat{\theta}_{0,\ldots,0}(h).$ Note that this local
approach allows to estimate successive derivatives of function $ f $. In this paper, only the estimation of $ f $ at a given point is studied. 

We note that similar locally parametric approach based on maximum likelihood estimators was recently proposed in
\cite{Katkovnik_Spokoiny08} and \cite{Polzehl_Spokoiny06} for {\it regular statistical models}.
But when the density of observations is discontinuous, the bayesian approach outperforms the maximum likelihood estimator. This phenomenon is well
known in parametric estimation (see \cite{Hasminskii_Ibragimov81}). Moreover, the establishing of statistical properties of bayesian
estimators requires typically much weaker assumptions than whose used for analysis of maximum likelihood  estimators.

\medskip

As we see our construction contains an extra-parameter $h$ to be
chosen. To make this choice we use quite standard arguments. First,
we note that in view of the definition of $ \bH_d({\beta},L,M) $ (below in Definition \ref{def_holder}), we have $ \forall f\in\bH_d({\beta},L,M), $
\begin{eqnarray*}
\exists\theta=\theta(f,y,h)\in [-M,M]^{D_b}\::\sup_{x\in
V_h(y)}\big|f(x)-f_\theta(x)\big|\leq Ldh^\beta.
\end{eqnarray*}
Remark that if $ f\in\bH_d({\beta},L,M) $, then $ \theta\in \Theta(A,M) $.
Thus, if $h$ is chosen sufficiently small, our original model (\ref{model}) is well approximated inside of $V_h(y)$ by the ``parametric" model
$$
\cY_i=f_\theta(X_i)\times U_i,\quad i=1,\ldots, nh^d,\quad nh^d\in\bN^*
$$
in which the {\it bayesian estimator} $\hat{\theta}$ is rate-optimal (See Theorem \ref{minimax}).

It is worth mentioning that the analysis of the deviation of $(X_i,\cY_i)_{i=1,...nh^d}$ from $Y^{(nh^d)}$  is not simple. Namely here requirements
$0<A\leq f(x)\leq M, \forall x\in[0,1]^{d},$ are used. This assumption, which seems not to be necessary, allows us to make the presentation of basic
ideas clear and to simplify routine computations (see also Remark \ref{rem_low_bound}).

Finally, $h=h_n({\beta},L)=(Ln)^{-1/(\beta+d)}$ is chosen as the solution of the following
minimization problem
\begin{equation}\label{equation minimiser}
Ldh^\beta+1/nh^d\rightarrow\min_{h}
\end{equation}
and we show that corresponding estimator
$\bar{f}^{h_n({\beta},L)}(y)$ is minimax for $f(y)$ on
$\bH_d({\beta},L,M,A)$ if $\beta\leq b$ (see Theorem \ref{minimax}).  Since the parameter $b>0$
can be chosen in arbitrary way, the proposed estimator is minimax for any given value of the parameter $\beta>0$. 

We remark that in regular statistical models, where linear methods are usually optimal, the choice of the bandwidth $ h $ is due to the relation
$$
Ldh^\beta+1/\sqrt{nh^d}\rightarrow\min_{h},
$$
with the solution $ h_L=(Ln)^{-1/(2\beta+d)} $.
This explains that the improvement of the rate of convergence, $ (1/n)^{\beta/(\beta+d)} $ compared to $ (1/n)^{\beta/(2\beta+d)} $, in the model with
the discontinuous density.


\medskip

\paragraph{Adaptive estimation.}The second part of the paper is devoted
to the adaptive minimax estimation over collection of isotropic
functional classes in the model (\ref{model}). At our knowledge, the
problem of adaptive estimation in the multiplicative regression with
the noise, having discontinuous density,  is not studied in the
literature.
\medskip

Well-known drawback of minimax approach is the dependence of the
minimax estimator on the parameters describing functional class on
which the maximal risk is determined.  In particular, the locally
bayesian estimator $\bar{f}^h(\cdot)$ depends obviously on the
parameters $A$ and $M$ via the solution of the minimization problem
(\ref{bayes-estimator}). Moreover $h_n({\beta},L)$ optimally chosen
in view of (\ref{equation minimiser}) depends explicitly on
${\beta}$ and $L$. To overcome this drawback the minimax adaptive
approach was proposed (see \cite{Lepski90}, \cite{Lepski91}, \cite{Lepski_Mammen_Spokoiny97}). The first
question arising in the adaptation (reduced to the problem at hand)
can be formulated as follows.

\smallskip

{\it Does there exist an estimator which would be  minimax on
$\bH({\beta},L,M,A)$ simultaneously for all values of ${\beta},L,A$
and $M$ belonging to some given subset of $\bR^{4}_{+}$ ?}

\smallskip

In section \ref{sectionAdaptive}, we show that the answer to this
question is \textsf{negative}, that is typical for the estimation of
the function at a given point (see \cite{Lepski_Spokoiny97}, \cite{Spokoiny96}, \cite{Tsybakov98}). This
answer can be reformulated in the following manner: the family of
rates of convergence $\big\{\varphi_n(\beta),\:
\beta\in\bR_+^*\big\}$ is \textsf{unattainable} for the problem
under consideration.

Thus, we need  to find another family of normalizations for maximal
risk which would be attainable  and, moreover, optimal in view of
some criterion of optimality. Nowadays, the most developed criterion
of optimality is due to \emph{Klutchnikoff} \cite{Klutchnikoff05}.

\medskip

We show that the family of normalizations, being optimal in view of
this criterion,
 is
\begin{equation}
\label{adaptive_rate}
\phi_{n}(\beta)=\left(\frac{\rho_{n}(\beta)}{n}\right)^{\frac{\beta}{\beta+d}},\qquad
\rho_{n}(\beta)=1+\ln\left(\frac{\varphi_{n}({\beta})}{\varphi_{n}({b})}\right),
\end{equation}
 whenever $\beta\in]0,b].$
 The factor $\rho_{n}$ can be considered as {\it price to
pay for adaptation} (see \cite{Lepski91}).
\medskip

The most important step in proving the optimality of the family
(\ref{adaptive_rate}) is to find an estimator, called {\it
adaptive},  which attains the optimal family of norma\-lizations.
 Obviously, we seek an estimator whose construction is {\it parameter-free}, i.e. independent of $\beta,L,A$ and $M$.
  In order to explain our estimation procedure let us make several remarks.

  First we note that the role of the constants $A,M$  and $\beta,L$ in the construction of the minimax estimator is quite different.
  Indeed, the constants $A,M$ are used in order to determine the set $\Theta\big(A,M\big)$ needed for the construction of the locally
  bayesian estimator, see (\ref{critere}) and (\ref{bayes-estimator}).
  However, this set does not depend on the localization parameter $h>0$, in other words, the quantities  $A$ and $M$ are
  not involved in the selection of optimal size of the local neighborhood given by (\ref{equation minimiser}).
  Contrary to that, the constants $\beta,L$ are used for the derivation of the optimal size of the local neighborhood
  (\ref{equation minimiser}), but they are not involved in the construction of the collection of locally bayesian
  estimators $\big\{\hat{f}^h, h>0\big\}.$

  Next remark explains how to replace the unknown quantities $A$ and $M$ in the definition of  $\Theta\big(A,M\big)$.
  Our first simple observation consists in the following: the estimator $\bar{f}^{h_n(\beta,L)}$ remains minimax
  if we replace $\Theta\big(A,M\big)$ in (\ref{critere}) and (\ref{bayes-estimator}) by
  $\Theta\big(\tilde{A},\tilde{M}\big)$ with any $0<\tilde{A}\leq A$ and $M\leq\tilde{M}<\infty$. It follows from obvious inclusion
  $\bH_d(\beta,L,A,M)\subseteq \bH_d(\beta,L,\tilde{A},\tilde{M}).$
  The next observation is less trivial and it follows from Proposition \ref{lemma1}. Put $h_{\max}=n^{-\frac{1}{b+d}}$ and define
  for any function $f$
\begin{equation}\label{borneInfSup}
  A(f)=\inf_{x\in V_{h_{\max}}(y)}f(x),\quad M(f)=\sum_{m=0}^{b}\sum_{p_1+\ldots+p_d=m}\left|\frac{\partial^{m}f(y)}{\partial x_1^{p_1}\cdots\partial
x_d^{p_d}}\right|.
\end{equation} The following agreement will be
used in the sequel:  if the function $f$ and  $m\geq 1$ be such that
$\partial^{m}f$  does not exist we will put formally
$\partial^{m}f=0$ in the definition of $M(f)$.

It remains to note that contrary to the
quantities $A$ and $M$ the functionals $A(f)$ and $M(f)$ can be
consistently estimated from the observation (\ref{model}) and let
$\hat{A}$ and $\hat{M}$ be  the corresponding estimators.
 The idea now is to determine  the collection of locally bayesian estimators  $\big\{\hat{f}^h, h>0\big\}$  by replacing $\Theta\big(A,M\big)$ in
 (\ref{critere}) and (\ref{bayes-estimator}) by the {\it random} parameter set $\hat{\Theta}$ which is defined as follows.
 $$
 \hat{\Theta}=\Theta\big(\hat A/2,4\hat M\big)=\left\{t\in\bR^{D_b}:\:\: 2t_{0,...,0}-\|t\|_1\geq 2^{-1}\hat{A},\:\: \|t\|_1\leq 4\hat{M}\right\}.
 $$
 In this context it is important to emphasize that the estimators  $\hat{A}$ and $\hat{M}$ are built from the same observation which is used for the
construction of the  family
 $\big\{\hat{f}^h, h>0\big\}$.

 \smallskip

 Contrary to all saying above, the constants $\beta$ and $L$ cannot be estimated consistently. In order to select an ``optimal" estimator from the
 family $\big\{\hat{f}^h, h>0\big\}$ we use  general adaptation
scheme due to \emph{Lepski} \cite{Lepski90}, \cite{Lepski92a}. To the best of our knowledge it is the first time when this method is applied in the
statistical model with multiplicative noise and discontinuous distribution. Moreover, except already mentioned papers
\cite{Katkovnik_Spokoiny08} and \cite{Polzehl_Spokoiny06}, Lepski's procedure is typically applied to the selection from the collection of linear
estimators (kernel estimators, locally polynomial estimator, etc.). In the present paper we apply this method to very complicated family of nonlinear
estimators, obtained by the use of bayesian approach on the random parameter set. It required, in particular, to establish the exponential inequality
for the deviation of locally bayesian estimator from the parameter to be estimated (Proposition \ref{lemma1}). It generalizes the  inequality
proved for the parametric model (see \cite{Hasminskii_Ibragimov81} Chapter 1, Section 5), this result seems to be new.

\paragraph{Simulations.} In the present paper we adopt the local parametric
approximation to a purely non parametric model. As it proved, this
strategy leads to the theoretically optimal statistical decisions.
But the minimax as well as the minimax adaptive approach are
asymptotical and it seems natural to check how proposed estimators
work for reasonable sample size. In the simulation study, we test
the bayesian estimator in the parametric and nonparametric cases. We
show that the \emph{adaptive} estimator approaches the \emph{oracle}
estimator. The \emph{oracle} estimator is selected from the family
$\left\{\hat f^{h},\,h>0\right\}$ under the hypothesis \emph{f} that
is known. We show that the bayesian estimator performs well starting
with $n\geq100$.

\bigskip

This paper is organized as follows. In Section \ref{sectionMinimax}
we present the results concerning minimax estimation and Section
\ref{sectionAdaptive} is devoted to the adaptive estimation. The
simulations are given in Section \ref{sectionSimu}. The proofs of
main results are proved in Section \ref{sectionUpper_bounds} (upper
bounds) and section \ref{sectionlower_bounds} (lower bounds).
Auxiliary lemmas are postponed to Appendix (Section \ref{appendixA}) contains the proofs of technical results.

\section{Minimax estimation on isotropic H\"older class}\label{sectionMinimax}In this section we present several results concerning minimax
estimation. First, we establish lower bound for minimax risk defined
on $\mathbb{H}_d({\beta},L,M,A)$ for any ${\beta},L,M$ and $A$. For
any $(p_1,...,p_d)\in\bN^d$ we denote $p=(p_1,...,p_d)$ and
$|p|=p_1+...+p_d$.

\begin{definition}
\label{def_holder} Fix $\beta>0$, $L>0$ and $M>0$ and  let
$\lfloor\beta\rfloor$ be the largest integer  strictly less than
$\beta$. The \emph{isotropic Hölder class} $\bH_d({\beta},L,M)$ is
the set of  functions $f:[0,1]^d\rightarrow\bR$ having on $[0,1]^d$
all partial derivatives of order $\lfloor\beta\rfloor$ and such that
$\forall x,y\in [0,1]^d$

\begin{eqnarray*}
\sum_{0\leq|p|\leq \lfloor\beta\rfloor}\sup_{x\in[0,1]^d}\left|\frac{\partial^{|p|}f(x)}{\partial
x_1^{p_1}\cdots\partial x_d^{p_d}}\right|&\leq& M,\\*[2mm]
\left|f(x)-\sum_{0\leq|p|\leq \lfloor\beta\rfloor} \frac{\partial^{|p|}f(y)}{\partial y_1^{p_1}\cdots\partial
y_d^{p_d}}\prod_{j=1}^d\frac{(x_j-y_j)^{p_j}}{p_j!}\right|&\leq&L\sum_{j=1}^d|x_j-y_j|^{\beta},
\end{eqnarray*}
where $x_j$ and $y_j$ are the $j$th components of $x$ and $y$.

\end{definition}
This definition implies that if $ f\in\mathbb{H}_d({\beta},L,M,A) $ (defined in the beginning of this paper), then $ A\leq A(f) $ and $ M(f)\leq M $,
where $ A(f) $ and $ M(f) $ are defined in (\ref{borneInfSup}).
\paragraph{Maximal and minimax risk on $\mathbb{H}_d({\beta},L,M,A)$.}To measure the performance of estimation procedures on
$\mathbb{H}_d({\beta},L,M,A)$ we will use minimax approach.

Let $\bE_f=\bE^{n}_f$ be the mathematical expectation with respect
to the probability law of the observation $Y^{(n)}$ satisfying
(\ref{model}). We define first the maximal risk on
$\mathbb{H}_d({\beta},L,M,A)$ corresponding to the estimation of the
function $f$ at a given point $y\in[0,1]^{d}$.

\noindent Let $\tilde{f}$ be an arbitrary estimator built from the
observation $Y^{(n)}$. Let $\forall q>0$
$$
R_{n,q}\big[\tilde{f},\mathbb{H}_d({\beta},L,M,A)\big]=\sup_{f\in\mathbb{H}_d({\beta},L,M,A)}\mathbb{E}_f\big|\tilde{f}(y)-f(y)\big|^q.
$$
The quantity
$R_{n,q}\big[\tilde{f},\mathbb{H}_d({\beta},L,M,A)\big]$ is called
{\it maximal risk} of the estimator $\tilde{f}$ on
$\mathbb{H}_d({\beta},L,M,A)$ and the {\it minimax risk} on
$\mathbb{H}_d({\beta},L,M,A)$ is defined as
$$
R_{n,q}\big[\mathbb{H}_d({\beta},L,M,A)\big]=\inf_{\tilde{f}}R_{n,q}\big[\tilde{f},\mathbb{H}_d({\beta},L,M,A)\big],
$$
where $\inf$  is taken over the set of  all  estimators.
\begin{definition}
\label{def_minimax} The normalizing sequence $\psi_n$ is called
minimax rate of convergence (MRT) and the estimator $\hat{f}$ is
called minimax (asymptotically minimax) if
\begin{eqnarray*}
\liminf_{n\to\infty}\psi^{-q}_n\:R_{n,q}\big[\hat{f},\mathbb{H}_d({\beta},L,M,A)\big]&>&0;\\
\limsup_{n\to\infty}\psi^{-q}_n\:R_{n,q}\big[\hat{f},\mathbb{H}_d({\beta},L,M,A)\big]&<&\infty.
\end{eqnarray*}
\end{definition}
\begin{theorem}\label{th:lower_bound}
\label{inf} For any ${\beta}>0$, $L>0$, $M>0,\:A>0$, $q\geq1$ and
$d\geq1$
\begin{eqnarray*}
\liminf_{n\to\infty}\varphi_{n}^{-q}(\beta)R_{n,q}\big[\mathbb{H}_d({\beta},L,M,A)\big]>0,\qquad
\varphi_{n}(\beta)=n^{-\frac{{\beta}}{{\beta}+d}}.
\end{eqnarray*}
\end{theorem}
\begin{remark}
\label{rem_low_bound} The obtained result shows that on
$\mathbb{H}_d\big({\beta},L,M,A\big)$  the minimax rate of
convergence cannot be faster than $n^{-\frac{{\beta}}{{\beta}+d}}$.
In view of the obvious inclusion
$\mathbb{H}_d({\beta},L,M,A)\subset\mathbb{H}_d({\beta},L,M)$ the
minimax rate of convergence on an isotropic Hölder class is also
bounded from below by $n^{-\frac{{\beta}}{{\beta}+d}}$.
\end{remark}
The next theorem shows how to construct the minimax estimator basing
on locally bayesian approach. Put $ \bar
h=(Ln)^{-\frac{1}{\beta+d}}$ and let $\bar{f}^{\bar
h}(y)=\hat{\theta}_{0,\ldots,0}\big(\bar h\big)$ is given by
(\ref{ensemble}), (\ref{critere}) and  (\ref{bayes-estimator}) with
$h=\bar{h}.$
\begin{theorem}\label{minimax}
Let  $\beta>0$, $L>0,\,M>0$ and $A>0$ be fixed. Then there exists
the constant $C_*$ such that for any $ n\in\bN^*$ satisfying
$n\bar h^d\geq\big(\lfloor\beta\rfloor+1\big)^d$
\begin{eqnarray*}
\varphi_{n}^{-q}(\beta)\:R_{n,q}\Big[\bar{f}^{\bar h}(y),\mathbb{H}_d({\beta},L,M,A)\Big]\leq C^*,\quad \forall q\geq1.
\end{eqnarray*}
The explicit form of $ C^* $ is given in the proof.
\end{theorem}
\begin{remark}
We deduce from Theorems \ref{inf} and \ref{minimax} that the
estimator $\bar{f}^{\bar h}(y)$ is minimax on
$\mathbb{H}_d({\beta},L,M,A)$.
\end{remark}

\section{Adaptive estimation on isotropic H\"older classes}\label{sectionAdaptive}This section is devoted to the adaptive estimation over the
collection of the classes
$\Big\{\mathbb{H}_d(\beta,L,M,A)\Big\}_{\beta,L,M,A}$. We will not
impose any restriction on  possible values of $L,M,A$, but we will
assume that $\beta\in (0,b]$, where $b$, as previously, is an
arbitrary a priori chosen integer.

We start with formulating the result showing that there is no
optimally adaptive estimator (here we follow the terminology
introduced in \cite{Lepski91}, \cite{Lepski92a}). It means that there
is no an estimator which would be minimax simultaneously for several
values of parameter $\beta$ even if all other parameters $L,M$ and
$A$ are supposed to be fixed. This result does not require any
restriction on $\beta$ as well.
\begin{theorem}
\label{th:nonexistence-adaptive-estimators} For any $\mathbb{B}\subseteq\bR^{+}\setminus \{0\}$ such that
$\text{card}(\mathbb{B})\geq 2$,
for any $\beta_1,\beta_2\in\mathbb{B}$ and any $L>0,\: M>0$, $A>0$
\begin{eqnarray*}
&&\liminf_{n\to\infty}\inf_{\tilde{f}}\Big[\varphi_{n}^{-q}(\beta_1)\:R_{n,q}\big(\tilde{f},\mathbb{H}_d(\beta_1,L,M,A)\big)\\
&&\qquad\qquad+
\varphi_{n}^{-q}(\beta_2)\:R_{n,q}\big(\tilde{f},\mathbb{H}_d({\beta}_2,L,M,A)\big)\Big]=+\infty,
\end{eqnarray*}
where $\inf$ is taken over all possible estimators.
\end{theorem}
 The assertion of Theorem \ref{th:nonexistence-adaptive-estimators} can be considerably specified if $\mathbb{B}=(0,b]$. To do that we will need the
following definition.
Let $ \Psi=\left\{\psi_n(\beta)\right\}_{\beta\in (0,b]} $ be a
given family of normalizations.
\begin{definition}\label{admissible}
The family $\Psi$ is called admissible  if there exist an estimator
$\hat{f}_n$ such that for \textsf{some} $L>0, M>0$ and $A>0$
\begin{equation}
\label{eq:def-admissibility} \limsup_{n\rightarrow\infty}
\psi^{-q}_n(\beta)\:R_{n,q}\big(\hat{f},\mathbb{H}_d(\beta,L,M,A)\big)<\infty,\:\:\forall
\beta\in (0,b].
\end{equation}
The estimator $\hat{f}_n$ satisfying (\ref{eq:def-admissibility}) is
called $\Psi$-\textsf{attainable}. The estimator $\hat{f}_n$ is
called $\Psi$-\textsf{adaptive} if (\ref{eq:def-admissibility})
holds  for \textsf{any} $L>0, M>0$ and $A>0$.
\end{definition}
Note that the result proved in  Theorem
\ref{th:nonexistence-adaptive-estimators} means that the family of
rates of convergence $\left\{\varphi_n(\beta)\right\}_{\beta\in
(0,b]}$ is not admissible. Denote by $\Phi$ the following family of
normalizations:
$$
\phi_{n}(\beta)=\left(\frac{\rho_{n}(\beta)}{n}\right)^{\frac{\beta}{\beta+d}},\quad
\rho_{n}(\beta)=1+\ln\left(\frac{\varphi_{n}(\beta)}{\varphi_{n}(b)}\right),\:\:\beta\in
(0,b].
$$
We remark that $\phi_n(b)=\varphi_n(b)$ and $\rho_{n}(\beta)\sim\ln
n$ for any $\beta\neq b$.
\begin{theorem}
\label{th:optimality-lower-bound} Let
$\Psi=\left\{\psi_n(\beta)\right\}_{\beta\in (0,b]}$ be an arbitrary
admissible family of normalizations.\\
\smallskip
{\bf I.} For any $\alpha\in (0,b]$ such that $\psi_n(\alpha)\neq\varphi_n(\alpha)$, there exists an admissible family
$\left\{\upsilon_n(\beta)\right\}_{\beta\in (0,b]}$ for which
$$
\lim_{n\to\infty} \upsilon_n(\alpha)\psi^{-1}_n(\alpha)=0.
$$
{\bf II.} If there exists $\gamma\in (0,b)$ such that
\begin{equation}
\label{eq:th-optimality-lower-bound} \lim_{n\to\infty}
\psi_n(\gamma)\phi^{-1}_n(\gamma)=0,
\end{equation}
then necessarily
\begin{eqnarray*}
&{\bf (a)}&\quad
 \lim_{n\to\infty} \psi_n(\beta)\phi^{-1}_n(\beta)>0, \quad
\forall \beta\in (0,\gamma);
\\*[2mm]
&{\bf (b)}&\quad \lim_{n\to\infty}
\left[\frac{\psi_n(\gamma)}{\phi_n(\gamma)}\right]\left[\frac{\phi_n(\beta)}{\psi_n(\beta)}\right]=0,\quad
\forall \beta\in (\gamma,b].
\end{eqnarray*}
\end{theorem}
Several remarks are in order.

We note that if the  family of normalizations $\Phi$ is admissible,
i.e. one can construct $\Phi$-attainable estimator, then $\Phi$ is
in an {\it optimal} family of normalizations in view of Kluchnikoff
criterion \cite{Klutchnikoff05}. It follows from the second assertion of the
theorem. We note however that a $\Phi$-attainable estimator may
depend on $L>0, M>0$ and $A>0$, and, therefore, this estimator  have
only theoretical interest. In the next section we construct
$\Phi$-adaptive estimator, which is, by its definition, fully
parameter-free. Moreover, this estimator obviously proves that
$\Phi$ is admissible, and, therefore, optimal as it was mentioned
above.

The assertions of Theorem \ref{th:optimality-lower-bound} allows us
to give rather simple interpretation of Kluchnikoff criterion.
Indeed, the first assertion, which is easily deduced from Theorem
\ref{th:nonexistence-adaptive-estimators}, shows that any admissible
family of norma\-lizations can be improved by another admissible
family  at any given point $\alpha\in (0,b]$ except maybe one. In
particular, it concerns the family $\Phi$ if it is admissible. On
the other hand, the second assertion of the theorem shows that there
is no  admissible family which would outperform the family $\Phi$ at
two points. Moreover, in view of $(\mathbf{b})$,
$\Phi$-adaptive (attainable) estimator, if exists, has the same
precision on $\mathbb{H}_d(\beta,L,M,A)$, $\beta<\gamma$, as any
$\Psi$-adaptive(attainable) estimator whenever $\Psi$ satisfies
(\ref{eq:th-optimality-lower-bound}). Additionally, $(\mathbf{a})$
implies that the gain in the precision provided by
$\Psi$-adaptive (attainable) estimator on
$\mathbb{H}_d(\gamma,L,M,A)$ leads automatically to much more losses
on $\mathbb{H}_d(\beta,L,M,A)$ for any  $\beta>\gamma$ with respect
to the precision provided by $\Phi$-adaptive(attainable) estimator.
We conclude that $\Phi$-adaptive(attainable) estimator outperforms
any $\Psi$-adaptive(attainable) estimator whenever $\Psi$ satisfies
(\ref{eq:th-optimality-lower-bound}). It remains to note that any
admissible family not satisfying
(\ref{eq:th-optimality-lower-bound}) is asymptotically equivalent to
$\Phi$.

\paragraph{Construction of $\Phi$-adaptive estimator.}As it was already mentioned in Introduction the construction of our estimation procedure
consists of several steps. First, we determine the set $\hat{\Theta}$, built from observation, which is used after that in order to define   the
family of locally bayesian estimators. Next, based on Lepski's method (see \cite{Lepski91} and \cite{Lepski_Mammen_Spokoiny97}), we propose
data-driven selection from this family.

{\it \underline{First step: Determination of parameter set.}} \: Put
$h_{\max}=n^{-\frac{1}{b+d}}$ and let $\tilde{\theta}$ be the
solution of the following minimization problem.
$$
\inf_{t\in\bR^{D_b}}\:\sum_{i:X_i\in V_{\max}(y)}^{n}\left[2Y_i-t\: K^\top\left(\frac{X_i-y}{h_{\max}}\right)\right]^{2},\quad
V_{\max}(y)=V_{h_{\max}}(y),
$$
where the $D_b$-dimensional vector $ K(z)=(z^{p}\::\:p\in\cP_b) $ and the sign $\top$
below means the transposition. Thus, $\tilde{\theta}$ is the local least squared estimator and its
explicit expression is given by
$$
\tilde{\theta}=2\left[\sum_{i:X_i\in V_{\max}(y)}^{n}K^\top\left(\frac{X_i-y}{h_{\max}}\right)
K\left(\frac{X_i-y}{h_{\max}}\right)\right]^{-1}\:\big[\cK_n(y)\big]^{\top}Y,
$$
where $Y=(Y_1,\ldots,Y_n)$ and
$\cK_n(y)=\left[K^\top\left(\frac{X_i-y}{h_{\max}}\right)\bI_{V_{\max}(y)}(X_i)\right]_{i=1,\ldots
n}$ is the design matrix. Put
$$
\displaystyle\tilde\delta_{p}=p_1!...p_d!\:h_{\max}^{-|p|}\:\tilde\theta_{p},\quad |p|\leq b.
$$
Introduce the following quantities
\begin{eqnarray}
\label{eq:def-estimators-for-A-M}
\hat{A}=\tilde{\delta}_{0\ldots,0},\qquad
\hat{M}=\big\|\tilde{\delta}\big\|_1,
\end{eqnarray}
and define the random parameter set as follows.
\begin{eqnarray}\label{eq:random-parameter-set}
\hat{\Theta}=\left\{t\in\bR^{D_b}:\:\: 2t_{0,...,0}-\|t\|_1\geq 2^{-1}\hat{A},\:\: \|t\|_1\leq 4\hat{M}\right\}.
\end{eqnarray}
{\it \underline{Second step: Collection of locally bayesian
estimators.}} \: Put
\begin{eqnarray}
 \label{critere-random}
    \hat{\pi}_h(t)&=&\int_{\hat{\Theta}}\|t-u\|_1\: L_h\big(u,Y^{(n)}\big)du;
 \\*[2mm]
  \label{bayes-estimator-random-set}
\hat{\theta}^*(h)&=&\arg\min_{t\in \hat{\Theta}}\hat{\pi}_h(t).
\end{eqnarray}
The family of locally bayesian estimator $\hat{\cF}$   is defined
now as  follows.
\begin{eqnarray}\label{eq:family-estimators}
\hat{\cF}=\left\{\hat{f}^h(y)=\hat{\theta}^*_{0,\ldots,0}(h),\: h\in
\big(0,h_{\max}\big]\right\}.
\end{eqnarray}
{\it \underline{Third  step: Data-driven selection from the
collection $\hat{\cF}$.}} \: Put
$$
h_k=2^{-k}h_{\max},\:\: k=0,\ldots, \mathrm{k}_n,
$$
where $\mathrm{k}_n$ is smallest integer such that $\displaystyle
h_{\mathrm{k}_n}\geq h_{\min}=\ln^{\frac{b}{d(b+d)}}n^{-1/d}$. Set
$$
\hat{\cF}^*=\left\{\hat{f}^{(k)}(y)=\hat{\theta}^*_{0,\ldots,0}(h_k),\:\:k=0,\ldots,
\mathrm{k}_n\right\}.
$$
We put $\hat{f}^*(y)=\hat f^{(\hat{k})}(y)$, where $\hat f^{(\hat{k})}(y)$ is
selected from $\hat{\cF}^*$ in accordance with the rule:
\begin{eqnarray}\label{indexe adaptive}
\hat k=\inf\left\{k=\overline{0,\mathrm{k}_n}:\:\:\big|\hat
f^{(k)}(y)-\hat f^{(l)}(y)\big|\leq \hat{M}S_n\big(l\big),\:\:
l=\overline{k+1,\mathrm{k}_n}\right\}.
\end{eqnarray}
Here we have used the following notations.
$$
S_n(l)=432D_b^3(32qd+16)\:\lambda^{-1}_n\big(h_l\big)\:\left[\frac{1+l\ln
2}{n\big(h_l\big)^d}\right],\quad l=0,1,\ldots,\mathrm{k}_n,
$$
and $\lambda_n(h)$ is the smallest eigenvalue of the matrix
\begin{eqnarray}\label{matrix_design}
\cM_{nh}(y)=\frac{1}{nh^{d}}\sum_{i=1}^{n}K^\top\left(\frac{X_i-y}{h}\right)
K\left(\frac{X_i-y}{h}\right)\bI_{V_{h}(y)}(X_i),
\end{eqnarray}
which is completely determined by the design points and by the number of observations. We will
prove that there exists a nonnegative real $\lambda$, such that $\lambda_n(h)\geq\lambda$ for any $n\geq 1$ and any $h\in \big[h_{\min},h_{\max}\big]$
(see Lemma \ref{matrixDP}).
\begin{theorem}
\label{lepski} Let an integer number $b>0$ be fixed. Then
 for any $\beta\in (0,b]$,  $L>0,\: M>0, A>0$ and $q\geq1$
$$
\limsup_{n\to\infty}\phi_{n}^{-q}(\beta)\:R_{n,q}\Big[\hat{f}^*(y),\mathbb{H}_d({\beta},L,M,A)\Big]<\infty.
$$
\end{theorem}
\begin{remark}
The assertion of the theorem means that the proposed estimator
$\hat{f}^*(y)$ is $\Phi$-adaptive. It implies in particular that the
family of normalizations $\Phi$ is admissible. This, together with
Theorem \ref{th:optimality-lower-bound} allows us to state the
optimality of $\Phi$ in view of Kluchnikoff criterion (see \cite{Klutchnikoff05}).
\end{remark}

\section{Simulation study}\label{sectionSimu}
\begin{figure}[!ht]
\center
\includegraphics[scale=0.3,width=\linewidth]{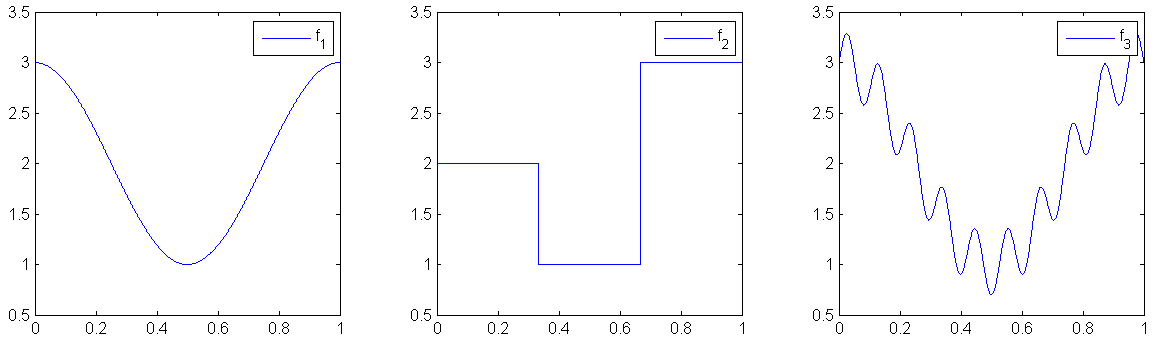} 
\caption{Test functions.}\label{test_function}
\end{figure}

We will consider the case  $d=1$.  The data are
simulated accordingly to the model (\ref{model}), where we use the
following functions (Figure \ref{test_function}).

Here $f_1(x)=\cos(2\pi x)+2$,
$f_2(x)=2.\bI_{[x\leq1/3]}+1.\bI_{[1/3<x\leq2/3]}+3.\bI_{[2/3<x]}$ and
$f_3(x)=\cos(2\pi x)+2+0.3\sin(19\pi x)$

\begin{figure}[!ht]
\center
\includegraphics[scale=0.3,width=\linewidth]{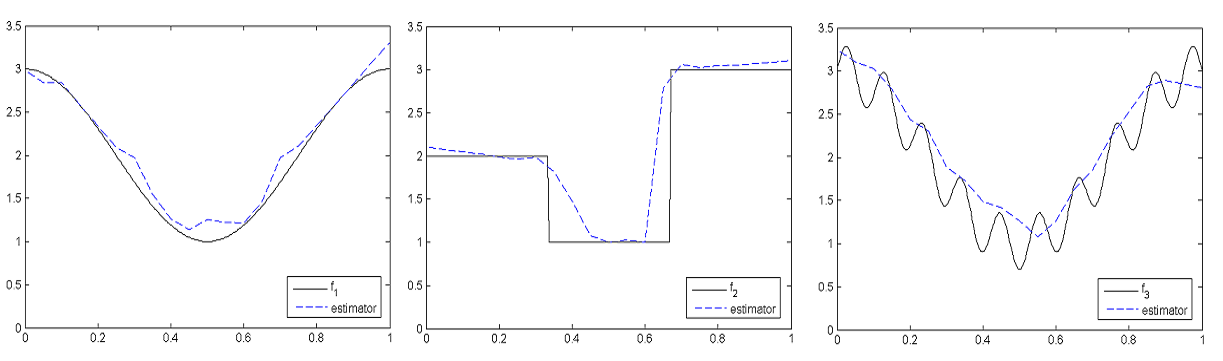} 
\caption{Examples of estimation with $ n=100 $.} \label{examples}
\end{figure}

To construct the  family of  estimators we use the linear
approximation ($b=2$), i.e. within the neighbourhoods of the given
size $h$, the locally bayesian estimator has the form
$$
\hat f^{h}(x)=\hat\theta_0+\hat\theta_1x,\quad x\in[0,1].
$$
We define the ideal (oracle) value of the parameter $\tilde
h=\tilde{h}(f)$ as the minimizer of the risk:
$$
\tilde h=\arg\inf_{h\in[1/n,1]}\mathbb{E}_f\big|\hat
f^{h}(y)-f(y)\big|.
$$
To compute it we apply  Monte-Carlo simulations (10000 repetitions).
Our first objective is to compare the risk provided by the "oracle"
estimator $\hat{f}^{\tilde{h}}(\cdot)$ and whose provided by the
adaptive estimator from Section \ref{sectionAdaptive}. Figure
\ref{examples} shows the deviation of the adaptive estimator from
the function to be estimated. In several points, for example in
$y=1/2$, we remark so-called over-smoothing phenomenon, inherent to
any adaptive estimator.

\paragraph{Oracle-adaptive ratio.}We compute the risks of the oracle and the adaptive estimator in 100
points of the interval $(0,1)$. The next tabular presents the mean
value of the ratio \textsf{oracle risk}/\textsf{adaptive risk}
calculated for the functions $f_1,f_2,f_3$ and $n=100,1000$.

\begin{figure}[!ht]
\begin{center}\label{tablo}
\begin{tabular}{|l||p{1.5cm}|p{2.5cm}||p{1.5cm}|p{2.5cm}|}
\hline \multicolumn{1}{|l|}{ } & \multicolumn{2}{|c|}{n = 100} & \multicolumn{2}{c|}{n = 1000} \\
\hline function & adaptive risk & oracle-adaptive ratio & adaptive risk & oracle-adaptive ratio\\
\hline $f_1$ & $0.13$ & $0.84$ & $0.03$ & $0.85$\\
\hline $f_2$ & $0.3$ & $0.71$ & $0.1$ & $0.75$\\
\hline $f_3$ & $0.28$ & $0.65$ & $0.2$ & $0.68$\\
\hline
\end{tabular}
\end{center}
\caption{Numeric values of risk.}
\end{figure}
\begin{figure}[!ht]
\center
\includegraphics[scale=0.3,width=\linewidth]{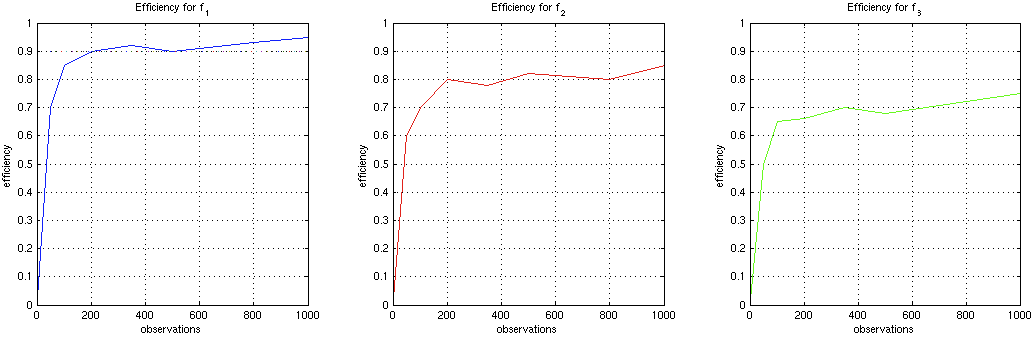} 
\caption{Efficiency of bayesian estimator for three test functions.} \label{ratio}
\end{figure}

Figure \ref{ratio} presents the "\textsf{oracle
risk}/\textsf{adaptive risk}" ratio as the function of the number of
observations $n$.

\begin{figure}[!ht]
\center
\includegraphics[scale=0.5]{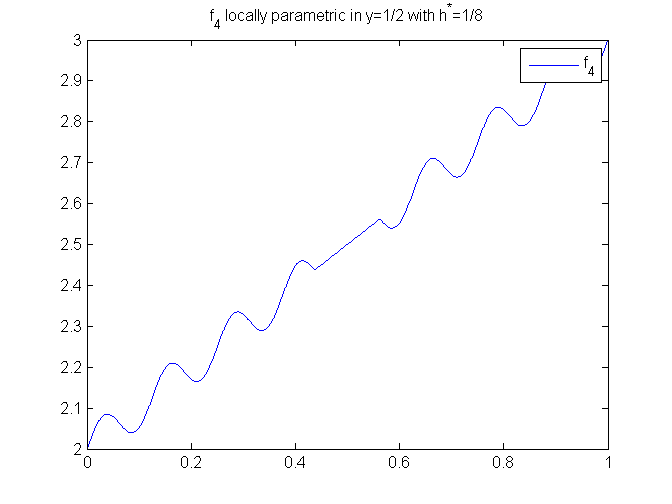} 
\caption{local parametric test function.}\label{f4}
\end{figure}

\paragraph{Adaptation versus parametric estimation.}We consider
the function $f_4$ (figure \ref{f4}), which is  linear inside the
neighborhood of size $h_*=1/8$ around point $1/2$ and simulate
$n=1000$ observations in accordance with the model (\ref{model}).
Using only the observations corresponding to the interval
$[3/8,5/8]$ we construct the bayesian estimator
$\hat{f}^{1/8}(1/2)$.

It is important to emphasize that this estimator is efficient \cite{Hasminskii_Ibragimov81} since the model is parametric. Our objective now is to
compare the risk of our adaptive estimator with the risk provided by the estimator $\hat{f}^{1/8}(1/2)$. We also try to understand how
far is  the localization parameter $h_{\hat{k}}$, inherent to the construction of our adaptive estimator, from the true value $1/8$.
We compute the risk of each estimator via Monte-Carlo method with $10000$ repetitions. For each repetition the procedure select the adaptive
bandwidth $h^{(j)}_{\hat k},\:j=1,...,10000$. We confirm once again the over-smoothing phenomenon since
$$
h^{(j)}_{\hat k}\sim0.1405>h_*=0.1250,\quad j=1,...,10000.
$$
Note however that the adaptive procedure selects the neighborhood of
the size which is quite close to the true one. We also compute the
risks of both estimators: ``\textsf{bayesian risk}"={\bf 0.0206} and
``\textsf{adaptive risk}"={\bf 0.0308}. We conclude that the
estimation accuracy provided by our adaptive procedure is quite
satisfactory.

\section{Proofs of main results: upper bounds}\label{sectionUpper_bounds} Let $\cH_n, n>1$
be the following subinterval of $(0,1)$.
\begin{equation}\label{Assumptionh}
\cH_n=\left[\frac{\big(b+1\big)\vee\big(\ln{n}\big)^{\frac{1}
{(d+d^{2})}}}{n^{1/d}},\left(\frac{1}{\ln{n}}\right)^{\frac{1}{b+d}}\right].
\end{equation}
Later on we will consider only the values of $h$ belonging to $\cH_n$. We start with establishing the exponential inequality for the deviation of
locally bayesian estimator $ \hat f^h(y) $ from $ f(y) $. The corresponding inequality is the basic technical to allowing to prove minimax and
minimax adaptive results.

\subsection{Exponential Inequality}\label{section_Auxi_result} Introduce the following notations. For any $h\in\cH_n$, put 
$\omega=\omega(f,y,h)=\big\{\omega_{p}:\:p\in\cP_b\big\}$, where $\omega_0=\omega_{0,...,0}=f(y)$ and
\begin{equation}\label{def_omega_taylor}
 \omega_{p}=\displaystyle \frac{\partial^{|p|}f(y)}{\partial y_1^{p_1}\cdots\partial y_d^{p_d}}\frac{h^{|p|}}{p_1!...p_d!},\:p\in\cP_b.
\end{equation}
Remind the agreement which we follow in the present paper: if the
function $f$ and vector $p$\: are such that $\partial^{|p|}f$  does not exist we put $\omega_{p}=0$.

Let $f_\omega(x)$, given by (\ref{polynome}), be the local
polynomial approximation of $f$ inside $V_h(y)$ and let $b_h$ be the
corresponding approximation error, i.e.
\begin{equation}\label{def_bias}
b_h=\sup_{x\in V_h(y)}\big|f_\omega(x)-f(x)\big|. 
\end{equation}
If $ f\in \bH_d(\beta,L,M),\:\beta>0 $, one could remark that $ b_h\leq Ldh^\beta $ by definition of $ \omega $ in (\ref{def_omega_taylor}) and $
\bH_d(\beta,L,M) $ in Definition \ref{def_minimax}.
Put also
\begin{equation}\label{borneexp}
\cN_h=b_h\times nh^d,\quad\cE(h)=\exp\left\{\frac{(1+6D_b^2)\cN_h}{6A(f)D_b^2}\right\}.
\end{equation}
Introduce the random events  $ G_{\hat M}=\big\{\big|\hat M-M(f)\big|\leq M(f)/2\big\}$ and
$G_{\hat A}=\big\{\big|\hat A-A(f)\big|\leq A(f)/2\big\}$ and put
$G=G_{\hat M}\cap G_{\hat A}$ where $ \hat A $ and $ \hat M $ are defined in (\ref{eq:def-estimators-for-A-M}), Section \ref{sectionAdaptive}.

Recall that $\lambda_n(h)$ (see Section \ref{sectionAdaptive}) is the smallest eigenvalue of the matrix
$$
\cM_{nh}(y)=\frac{1}{nh^{d}}\sum_{i=1}^{n}K^\top\left(\frac{X_i-y}{h}\right)
K\left(\frac{X_i-y}{h}\right)\bI_{V_{h}(y)}(X_i),
$$
and $ K(z) $ is the $D_b$-dimensional vector of the monomials $ z^p,\:p\in\cP_b $.

\begin{proposition}\label{lemma1}
For any $h\in\cH_n$ and any $f$ such that
$A(f)>A$ and $M(f)<M$, then $\forall\varepsilon>144MD_b(1\vee\cN_h)/A\lambda_n(h)$
\begin{eqnarray*}
\mathbb{P}_f\left(nh^d\big|\hat{f}^h(y)-f(y)\big|\geq \varepsilon,\:G\right) \leq
\mB\big(A(f),M(f)\big)\cE(h)\exp\left\{-\frac{\lambda_n(h)\:\varepsilon}{432M(f)\:D_b^3}\right\},
\end{eqnarray*}
where $\hat{f}^h(y)\in\hat\cF$ defined in
(\ref{eq:family-estimators}). The explicit expression of the
function $\mB(\cdot,\cdot)$ is given in the beginning of the proof
of the proposition.
\end{proposition}
The next proposition provides us with upper bound for the risk of a
locally bayesian estimator.
\begin{proposition}\label{2}
For any $n\in\bN^*$, $h\in\cH_n$ and any $f\in\bH_d(\beta,L,M,A)$,
then $\exists\lambda>0$ such that $ \lambda_n(h)\geq\lambda $ and
\begin{eqnarray*}
\bE_f\big|\hat{f}^h(y)-f(y)\big|^q\bI_G \leq C^*_q\big(A(f),M(f)\big)\left[\frac{1\vee Ld\:nh^{\beta+d}}{nh^d}\right]^{q},\quad
q\geq1,
\end{eqnarray*}
where
$$
C^*_q(a,m)=\frac{1}{q}\left[\frac{432m D_b^3(1+6D_b^2)}{3\lambda aD_b^2}\right]^q+\left[864m\lambda^{-1} D_b^3\right]^q\mB(a,m)\Gamma(q),\quad a,m>0,
$$
$ \Gamma(\cdot) $ is the well-known Gamma function.
\end{proposition}
\begin{remark}\label{remark_set_norandom}
The analysis of the proof of Proposition \ref{lemma1} allows to assert the following inequality
$$
\mathbb{P}_f\left(nh^d\big|\bar{f}^h(y)-f(y)\big|\geq \varepsilon\right) \leq
\mB\big(A,M\big)\cE(h)\exp\left\{-\frac{\lambda_n(h)\:\varepsilon}{432M\:D_b^3}\right\},
$$
where $ \bar f^h(y) $ is locally bayesian estimator which is the minimizer in (\ref{critere}). 

Thus, the latter inequality can be viewed an analogue of the result of Proposition \ref{lemma1} when $ A $ and $ M $ are known. By the same reasons,
we have 
\begin{eqnarray*}
\bE_f\big|\bar{f}^h(y)-f(y)\big|^q \leq C^*_q(A,M)\left[\frac{1\vee Ld\:nh^{\beta+d}}{nh^d}\right]^{q},\quad
q\geq1.
\end{eqnarray*}
\end{remark}

\subsection{Proof of Proposition \ref{lemma1}}
Before to start with the proof, let us breafly discuss its ingredients.

\paragraph{Discussion.}\hspace{0.2cm}

{\bf I.} First, the obvious inclusion (remind that $\hat\theta^*(h)$ minimizes $\hat\pi_h$ defined in (\ref{critere-random}))
$$
\left\{nh^d\big\|\hat\theta^*(h)-\theta\big\|_1\geq\varepsilon\right\}\subseteq
\left\{\inf_{nh^d\|t-\theta\|_1\geq\varepsilon}\hat{\pi}_h(t)\leq\hat{\pi}_h(\theta)\right\}.
$$
allows us to reduce the study of the deviation of $ \hat\theta^*(h)$ from $ \theta $ to the study of the behaviour of $ \hat{\pi}_h $.

{\bf II.} We note that $ \hat{\pi}_h $ is the integral functional of the pseudo-likelihood $ L_h $. As the consequence, the behaviour of $ \hat{\pi}_h
$ is completely determined by this process. Following \cite{Hasminskii_Ibragimov81} (Chapter 1, Section
5, Theorem 5.2), where similar problems were studied under parametric model assumption, we introduce the stochastic process 
$$
Z_{h,\theta}(u)=\frac{L_h\big(\theta+(nh^d)^{-1}u,Y^{(n)}\big)}{L_h\big(\theta,Y^{(n)}\big)}.
$$
defined on $ \Upsilon_n=\left\{u\in\bR^{D_b}\::\:u=nh^d(t-\theta),\:t\in\Theta(A(f)/4,9M(f))\right\}.$

Here, the vector $\theta=\theta(f,y,h)=\big\{\theta_{p}:\:p\in\cP_b\big\}$ is defined as follows.
$$
\theta_0=\theta_{0,...,0}=\omega_0+b_h,\quad \theta_{p}=\omega_{p},\:\:|p|\neq 0,
$$
where $ \omega $ is the coefficients of Taylor polynomial defined in (\ref{def_omega_taylor}). The definition of $b_h$ implies obviously
\begin{eqnarray}\label{controlbias}
f_\theta(x)&\geq& f(x),\quad \forall x\in V_h(y).
\end{eqnarray}
As it was noted in \cite{Hasminskii_Ibragimov81} (Chapter 1, Section 5, Theorem 5.2) the following properties of the process $ Z_{h,\theta} $ are
essential for the study of $ \hat{\pi}_h $:
\begin{itemize}
 \item Hölder continuity of its trajectories;
 \item the rate of its decay at infinity.
\end{itemize}
The exact statements are formulated in Lemma \ref{propertyZ} below.

{\bf III.} As it was shown in \cite{Hasminskii_Ibragimov81} (Chapter 1, Section 5, Theorem 5.2) in parametric situation the mentioned above properties
$ Z_{h,\theta} $ provide with the desirable properties of the process 
$$
z_h(u)=\frac{Z_{h,\theta}(u)}{\int_{\hat
\Upsilon_n}Z_{h,\theta}(v)dv},\quad u\in \hat{\Upsilon}_n:=nh^d\big(\hat\Theta-\theta\big),
$$
where the set $ \hat\Theta $ is defined in (\ref{eq:random-parameter-set}). The exact statements are given in Assertions 1 and 2. The latter process
is important in view of the following inclusion 
$$
\left\{nh^d\big|\hat{f}^h(y)-f(y)\big|\geq\varepsilon\right\}\subseteq\left\{\int_{\hat{\Upsilon}_n(r)}\|u\|_1z_{h}(u)du>\frac
{r}{2}\right\}.
$$

\paragraph{Auxiliary Lemma.}
First, we note that in view of (\ref{controlbias}), the event $ Y_i\leq f_{\theta}(X_i) $ is always realized, because $ Y_i\leq f(X_i)\leq
f_\theta(X_i) $. Hence, $ Z_{h,\theta} $ can be rewritten
\begin{equation}\label{not_processZ}
 Z_{h,\theta}(u)=\prod_{i:\:X_i\in
V_h(y)}{\frac{f_\theta(X_i)}{f_{\theta+u(nh^d)^{-1}}(X_i)}}\:\bI_{\left[Y_i\leq\:f_{\theta+u(nh^d)^{-1}}(X_i)\right]},
\quad u\in \Upsilon_n.
\end{equation}

\begin{lemma}\label{propertyZ}
For any $f\in \bH_d(\beta,L,M,A)$ and $h\in\cH_n$

\begin{enumerate}
 \item $ \displaystyle\sup_{u_1,u_2\in \Upsilon_n}\|u_1-u_2\|_1^{-1}\bE_f\big|Z_{h,\theta}(u_1)-Z_{h,\theta}(u_2)\big|\leq\cC_h,$
\item $ \displaystyle\bE_f Z^{1/2}_{h,\theta}(u)\leq e^{-g_h\big(||u||_1\big)},\:\:\forall u\in \Upsilon_n,$
\item $ \displaystyle\mathbb P_f\left\{\int_{[0,\delta]^{D_b}} Z_{h,\theta}(u)du<\frac{\delta^{D_b}}{2}\right\}<2{\cC_h\delta},\:\:\forall
\delta>0.$
\end{enumerate}
where
$$
\cC_h=8\big(1\vee
D_bA^{-1}(f)\big)\exp\left\{1+\cN_h/A(f)\right\},\quad
g_h(a)=\frac{\lambda_n(h)a}{18M(f)D_b}-\frac{\cN_h}{A(f)},
$$
with $ a>0 $ and $\lambda_n(h)$ is the smallest eigenvalue of the matrix $ \cM_{nh}(y) $ defined in (\ref{matrix_design}).
\end{lemma}

\paragraph{Proof of  Proposition \ref{lemma1}.}
Define for any $ u>0$ and $v>0$
\begin{eqnarray}\label{bornepui}
\mB(a,m)=\sup_{z\geq 0}16e\big(1\vee D_b
a^{-1}\big)\Sigma(m)\left[\cB_z+6\right]\exp\left\{-\frac{\lambda z}{432 v
D_b^3}\right\},
\end{eqnarray}
where
$\cB_z=z^{D_b+1}+2(2z+2\big)^{2D_b}+5+D_b\left(z^{D_b}+(2z+2\big)^{\frac{D_b}{2}-1}\right)$, $\lambda>0$ is defined such that: $
\lambda_n(h)\geq\lambda $ for any $n\in\bN^*$, $h\in\cH_n$ (for more details, see Lemma \ref{matrixDP}) and
$$
\Sigma(v)=\frac{c^2(v)\big(3-c(v)\big)}{\big(1-c(v)\big)^3},\quad
c(v)=\exp\big\{-(54 vD_b^2)^{-1}\big\},\quad v>0.
$$
\begin{assertion}\label{Asser_P1.1}
 For any $\varepsilon>0$, and for all $ r $ such that $ 0<r<\varepsilon/3 $, we assume
$$
\mathbb
P_f\left(nh^d\big|\hat{f}^h(y)-f(y)\big|\geq\varepsilon,G\right)\leq2\bP_f\left(\int_{\hat{\Upsilon}_n(r)}\|u\|_1z_{h}(u)du>\frac
{r}{2},G\right).
$$
\end{assertion}

\begin{assertion}\label{Asser_P1.2}
 For all $h\in\cH_n$ and any $f$ such that $A(f)>A$ and
$M(f)<M$, then for any $a>32MD_b(1\vee\cN_h)/(\lambda A)$
\begin{eqnarray*} 
\mathbb E_f\left[\int_{\hat \Upsilon_n\cap\big\{\|u\|_1>a\big\}} \|u\|_1z_{h}(u)\:du\:\bI_G\right]
\leq a\Sigma\big(M(f)\big)\: \cB_a\:\cC_h\exp\left\{-\frac{1}{6D_b}g_h(a)\right\},
\end{eqnarray*}
where $ g_h(\cdot) $ is defined in Lemma \ref{propertyZ}.
\end{assertion}

$1^0$. Suppose that Assertions 1 and 2 are proved.
Then, in view of Assertion \ref{Asser_P1.2}, choosing $r=\varepsilon/4$ we get
$$
\mathbb E_f\int_{\hat
\Upsilon_n\cap(\|u\|_1>\varepsilon/4)}\|u\|_1z_{h}(u)\:\bI_G\:du\leq
\frac{\varepsilon}{4}\:\Sigma\big(M(f)\big)\:\cB_{\varepsilon/4}\:\cC_h\:e^{-\frac{1}{6D_b^2}g_h(\varepsilon/4)}.
$$
Using the Tchebychev inequality, we have in view of the last inequality
\begin{eqnarray*}
\bP_f\left(\int_{\hat \Upsilon_n\cap(\|u\|_1>\varepsilon/4)}\|u\|_1z_{h}(u)du>\frac{\varepsilon}{8},\:G\right)
\leq 2\Sigma\big(M(f)\big)\:\cB_{\varepsilon/4}\:\cC_h\:e^{-\frac{1}{6D_b^2}g_h(\varepsilon/4)}.
\end{eqnarray*}
The assertion of Proposition \ref{lemma1} follows now from the last inequality, Assertion \ref{Asser_P1.1} and the definitions of $\cC_h, g_h(\cdot)$
and the function $\mB(\cdot,\cdot)$.

\bigskip

$2^0$. Now, we will prove Assertion \ref{Asser_P1.1}.
The definition of $\hat\theta^*(h)$ and $\theta=\theta(f,y,h)$
implies  $\forall\varepsilon>0$
\begin{eqnarray}\label{Cproba}
\mathbb{P}_f\left(nh^d\big|\hat{f}^h(y)-f(y)\big|\geq\varepsilon,G\right)&\leq&\mathbb{P}_f\left(nh^d\big|\hat\theta_0^*(h)-\theta_0\big|\geq
\varepsilon,G\right)\nonumber\\
&\leq&\mathbb{P}_f\left(nh^d\big\|\hat\theta^*(h)-\theta\big\|_1\geq
\varepsilon,G\right).
\end{eqnarray}
Some remarks are in order. First, it is easily seen that
$\theta\in\Theta\big(A(f),3M(f)\big)$. Therefore, if the event $G$
holds then $ \theta\in\hat{\Theta}. $ Remind  also that
$\hat\theta^*(h)$ minimizes $\hat\pi_h$  defined in (\ref{critere-random})
and, therefore, the following inclusion holds since
$\hat\theta^*(h)\in\hat{\Theta}.$
\begin{eqnarray}\label{eq1:proof-prop1}
\left\{\left(nh^d\big\|\hat\theta^*(h)-\theta\big\|_1\geq\varepsilon\right)\cap
G\right\}\subseteq
\left\{\left(\inf_{nh^d\|t-\theta\|_1\geq\varepsilon}\hat{\pi}_h(t)\leq
\hat{\pi}_h(\theta)\right)\cap G\right\}.
\end{eqnarray}
Moreover,
\begin{eqnarray*}
\hat\pi_h(t)&=&(nh^d)^{-1}\int_{\hat\Theta}\left\|nh^d(t-u)\right\|_1L_h\big(u,Y^{(n)}\big)du\\
&=&(nh^d)^{-D_b-1}\int_{\hat \Upsilon_n}\left\|nh^d(t-\theta)-u\right\|_1L_h\big(\theta+u(nh^d)^{-1},Y^{(n)}\big)du\\
&=&(nh^d)^{-D_b-1}L_h\big(\theta,Y^{(n)}\big)\int_{\hat
\Upsilon_n}\left\|nh^d(t-\theta)-u\right\|_1Z_{h,\theta}(u)du.
\end{eqnarray*}
Hence, $\tau_n=nh^d\big(\hat\theta^*(h)-\theta\big)$ is the minimizer
of
$$
\chi_n(s)=\int_{\hat
\Upsilon_n}\big\|s-u\big\|_1\frac{Z_{h,\theta}(u)}{\int_{\hat
\Upsilon_n}Z_{h,\theta}(v)dv}du
$$
and we obtain from (\ref{Cproba}) and (\ref{eq1:proof-prop1}) for
any $\varepsilon>0$
\begin{eqnarray}\label{eq2:proof-prop1} 
\mathbb P_f\left(\big\|nh^d(\hat\theta^*(h)-\theta)\big\|_1>\varepsilon,G\right)
\leq\mathbb P_f\left(\inf_{\|s\|_1>\varepsilon}\chi_n(s)\leq\chi_n(0),G\right).
\end{eqnarray}
 Let $0<r<\varepsilon/3,$ be a
number whose choice will be done later. We have
$$
\chi_n(0)\leq r\int_{\hat \Upsilon_n\cap(\|u\|_1\leq
r)}z_{h}(u)du+\int_{\hat \Upsilon_n\cap(\|u\|_1>r)}\|u\|_1z_{h}(u)du.
$$
Note also that
\begin{eqnarray*}
\inf_{\|s\|_1>\varepsilon}\chi_n(s)&\geq&\inf_{\|s\|_1>\varepsilon}\left[\int_{\hat
\Upsilon_n\cap(\|u\|_1\leq
r)}\big(\|s\|_1-\|u\|_1\big)z_h(u)du\right]\\
&\geq&(\varepsilon-r)\int_{\hat \Upsilon_n\cap(\|u\|_1\leq r)}z_h(u)du.
\end{eqnarray*}
It yields in particular
\begin{eqnarray*}
&&\chi_n(0)-\inf_{\|s\|_1>\varepsilon}\chi_n(s)\\
&&\quad\leq-(\varepsilon-2r)\int_{\hat \Upsilon_n\cap(\|u\|_1\leq
r)}z_{h}(u)du
+\int_{\hat \Upsilon_n\cap(\|u\|_1>r)}\|u\|_1z_{h}(u)du.
\end{eqnarray*}
Thus, $\forall r\in(0,\varepsilon/3)$
\begin{eqnarray}\label{borne_2proba}
&&\mathbb P_f\left(\chi_n(0)-\inf_{\|s\|_1>\varepsilon}\chi_n(s)>0,G\right)\nonumber\\
&&\leq\bP_f\left(\int_{\hat
\Upsilon_n\cap(\|u\|_1>r)}\|u\|_1z_{h}(u)du>(\varepsilon-2r)\int_{\hat
\Upsilon_n\cap(\|u\|_1\leq
r)}z_{h}(u)du,G\right)\nonumber\\
&&\leq\bP_f\left(\int_{\hat
\Upsilon_n\cap(\|u\|_1>r)}\|u\|_1z_{h}(u)du>r/2,G\right)\nonumber\\
&&\quad+\mathbb P_f\left((\varepsilon-2r)\int_{\hat
\Upsilon_n\cap(\|u\|_1\leq r)}z_{h}(u)du<r/2,G\right).
\end{eqnarray}
We note that the second term in (\ref{borne_2proba}) can be control by the first one whenever $0<r<\varepsilon/3$.  Indeed,
putting $\hat{\Upsilon}_n(r) =\hat \Upsilon_n\cap(u\in\bR^{D_b}\::\:\|u\|_1> r)$ we get
\begin{eqnarray*}
&&\mathbb P_f\left((\varepsilon-2r)\int_{\hat \Upsilon_n\cap(\|u\|_1\leq
r)}z_{h}(u)du<r/2,G\right)\\
&&\leq\mathbb P_f\left(r\int_{\hat
\Upsilon_n}Z_{h,\theta}(v)dv-r\int_{\hat{\Upsilon}_n(r)}Z_{h,\theta}(u)du<\frac{r}{2}\int_{\hat
\Upsilon_n}Z_{h,\theta}(v)dv,G\right)\\
&&\leq\mathbb
P_f\left(r\int_{\hat{\Upsilon}_n(r)}Z_{h,\theta}(u)du>\frac{r}{2}\int_{\hat
\Upsilon_n}Z_{h,\theta}(v)dv,G\right)\\
&&\leq\mathbb
P_f\left(\int_{\hat{\Upsilon}_n(r)}\|u\|_1z_{h}(u)du>r/2,G\right).
\end{eqnarray*}
The last inequality together with (\ref{Cproba}), (\ref{eq2:proof-prop1}) and (\ref{borne_2proba}) yields
$$
\mathbb
P_f\left(nh^d\big|\hat{f}^h(y)-f(y)\big|\geq\varepsilon,G\right)\leq2\bP_f\left(\int_{\hat{\Upsilon}_n(r)}\|u\|_1z_{h}(u)du>\frac
{r}{2},G\right).
$$

\smallskip

$3^0$. Now, let us prove Assertion \ref{Asser_P1.2}.
Put $\Upsilon_n(a)=\Upsilon_n\cap(u\in\bR^{D_b}\::\:\|u\|_1> a)$ for all $ a>0 $ and $\Omega_\upsilon=\Upsilon_n(\upsilon)\setminus
\Upsilon_n(\upsilon+1)$ for any $\upsilon\geq a$.
Introduce the following notations.
$$
\cI_\upsilon=\int_{\Omega_\upsilon} Z_{h,\theta}(u)du,\qquad \cQ_\upsilon=\frac{\int_{\hat \Upsilon_n\cap\:
\Omega_\upsilon}Z_{h,\theta}(u)du}{\int_{\hat
\Upsilon_n}Z_{h,\theta}(u)du}.
$$
Fix  $ T>0 $ whose choice will be done later.
Consider the minimal number $ N(\Omega_{\upsilon},1/T) $ of balls of radius $ 1/T $ that are needed to cover the set $ \Omega_{\upsilon} $. Denote
$u^j$ is the center of each ball. Since $ \Omega_{\upsilon} $ is a compact of $ \bR^{D_b} $, it implies  $N(\Omega_{\upsilon},1/T)\leq
(v+1)^{D_b}T^{D_b} $. Introduce the non-intersecting parts $
\Delta_1,\Delta_2,\Delta_3,\ldots $ as follows: $\Delta_1=\left\{u\in\Omega_{\upsilon}\::\: \|u-u^1\|_1\leq 1/T\right\}$ and
$$
\Delta_j=\left\{u\in\Omega_{\upsilon}\::\: \|u-u^j\|_1\leq
1/T\right\}\setminus\bigcup_{i=1}^{j-1}\Delta_j,\quad j=2,\ldots,N(\Omega_{\upsilon},1/T).
$$
Put $S_\upsilon=\sum_j\int_{\Delta_j}Z_{h,\theta}(u^j)du$ and note that $ S_\upsilon $ is stepwise approximation of $ \cI_\upsilon $.

\bigskip

{\it{ Control of $ \cI_\upsilon $.}} Remind that $ \Omega_{\upsilon}=\bigcup_{j=1}^{N(\Omega_{\upsilon},1/T)}\Delta_j $ and denote by
$\displaystyle
|\Omega_\upsilon| $ the volume of $ \Omega_\upsilon $. We get
for any $\sigma>0$
\begin{eqnarray*}
\mathbb P_f\big(S_\upsilon>\sigma\big)&\leq&\mathbb
P_f\left(\max_jZ_{h,\theta}^{1/2}(u^j)\sqrt{|\Omega_\upsilon|}>\sqrt\sigma\right)\\
&\leq& \sum_j \mathbb
P_f\left(Z_{h,\theta}^{1/2}(u^j)>\sqrt{|\Omega_\upsilon|}\sqrt\sigma\right).
\end{eqnarray*}
Note that the number of summands on the right-hand side of the last
inequality does not exceed $(v+1)^{D_b}T^{D_b}$. Applying Tchebychev inequality and
Lemma \ref{propertyZ} ({\it 2}), we obtain
\begin{equation}\label{5.11}
\mathbb P_f\big(S_\upsilon>\sigma\big)\leq
(v+1)^{D_b}T^{D_b}\sqrt{|\Omega_\upsilon|}\sigma^{-1/2}
e^{-g_h(\upsilon)}.
\end{equation}
In view of to Lemma \ref{propertyZ} ({\it 1}),
\begin{eqnarray*}
\bE_f\big|S_\upsilon-\cI_\upsilon\big|\leq\sum_j\int_{\Delta_j}\mathbb
E_f\big|Z_{h,\theta}(u)-Z_{h,\theta}(u^j)\big|du\leq{\cC_h}\sum_j\int_{\Delta_j}\|u-u^j\|_1du.
\end{eqnarray*}
By definition of $ \Delta_j $, each summand does not exceed $ \int_{\Delta_j}T^{-1}du$, therefore,
\begin{equation}\label{5.12}
\mathbb E_f\big|S_\upsilon-\cI_\upsilon\big|\leq\cC_h\sqrt{|\Omega_\upsilon|}T^{-1}.
\end{equation}
One has
\begin{eqnarray*}
\mathbb P_f\big(\cI_\upsilon>2\sigma\big)\leq \mathbb
P_f\big(S_\upsilon>\sigma\big)+\mathbb
P_f\big(\big|S_\upsilon-\cI_\upsilon\big|>\sigma\big).
\end{eqnarray*}
Using (\ref{5.11}), (\ref{5.12}) and applying  Tchebychev
inequality, we get
\begin{eqnarray}\label{5.10}
\mathbb P_f\big(\cI_\upsilon>2\sigma\big)\leq (v+1)^{D_b}T^{D_b}\sqrt{|\Omega_\upsilon|}\sigma^{-1/2}
 e^{-g_h(\upsilon)}+\cC_h\sqrt{|\Omega_\upsilon|}T^{-1}\sigma^{-1}.
\end{eqnarray}

\medskip

{\it Control of $ \cQ_\upsilon $.} Set $ \bA = \left\{\int_{\hat \Upsilon_n}Z_{h,\theta}(u)du<{\delta^{D_b}}/{2}\right\} $. Since $
\cQ_\upsilon\leq1 $ we obtain for any $\delta>0,\:\sigma>0$
\begin{eqnarray*}
\bE_f\cQ_\upsilon&=&\bE_f\left[\cQ_\upsilon\:\bI_{\bA}+\cQ_\upsilon\:\bI_{\cI_\upsilon>2\sigma,\bA^c}+\cQ_\upsilon\bI_{
\:\cI_\upsilon\leq2\sigma, \bA^c }\right]\bI_G \\
&\leq&\mathbb P_f\left(\bA,G\right)+\mathbb
P_f\big(\cI_\upsilon>2\sigma\big)+4\delta^{-D_b}\sigma.
\end{eqnarray*}
Under the event $ G $, remark that $ [0,\delta]^{D_b}\subseteq nh^d\big(\Theta\big(A(f),2M(f)\big)-\theta\big)\subseteq\hat \Upsilon_n $ for any $
\delta\leq \big(2M(f)-A(f)\big) $. Using to Lemma \ref{propertyZ} ({\it 3}) and the inequality (\ref{5.10}), we have
\begin{eqnarray*}
\mathbb
E_f\cQ_\upsilon\leq2{\cC_h\delta}+T^{D_b}\sqrt{|\Omega_\upsilon|}\sigma^{-1/2}
e^{-g_h(\upsilon)}+\cC_h\sqrt{|\Omega_\upsilon|}T^{-1}\sigma^{-1}+4\delta^{-D_b}\sigma.
\end{eqnarray*}
Choosing $T=\exp\left\{\frac{1}{2D_b}g_h(\upsilon)\right\}$, $\sigma=\exp\left\{-\frac{1}{3D_b}g_h(\upsilon)\right\}$ and
$\delta=\exp\left\{-\frac{1}{6D_b^2}g_h(\upsilon)\right\}$, we obtain
$$
\mathbb
E_f\cQ_\upsilon\leq\left[2{\cC_h}+\sqrt{|\Omega_\upsilon|}\left((v+1)^{D_b}+\cC_h\right)+4\right]\exp\left\{-\frac{1}{6D_b^2}g_h(\upsilon)\right\}.
$$

\medskip

{\it{Conclusion of the proof of Assertion \ref{Asser_P1.2}.}} Simplest algebra shows that $
\sqrt{|\Omega_\upsilon| }\leq\big(2\upsilon+2\big)^{\frac{D_b}{2}} $, we get
\begin{eqnarray}\label{Qtilde}
\mathbb E_f\cQ_\upsilon&\leq&
\big[\upsilon^{D_b+1}+2(2\upsilon+2\big)^{2D_b}+5\big]\:\cC_h\:\exp\left\{-\frac{1}{6D_b^2}g_h(\upsilon)\right\},
\end{eqnarray}
Note that if the event $G$ is realized then $\hat \Upsilon_n(a)\subseteq \Upsilon_n(a)=\bigcup_{j=0}^\infty\Omega_{a+j}. $ we
obtain in view of (\ref{Qtilde})
\begin{eqnarray*}
\bE_f\int_{\hat \Upsilon_n\cap(\|u\|_1>a)}\|u\|_1z_{h}(u)\:\bI_G\:du&\leq&\sum_{j=0}^\infty\big(a+j+1\big)\bE_f\cQ_{a+j}\\
&=&\Sigma\big(M(f)\big)\,a\,\cB_a\,\cC_h\,\exp\left\{-\frac{1}{6D_b^2}g_h(a)\right\}.
\end{eqnarray*}
where we have put
$\cB_a=a^{D_b+1}+2(2\upsilon+2\big)^{2D_b}+5+D_b\left(a^{D_b}+(2a+2\big)^{\frac{D_b}{2}-1}\right)$.
\epr
\subsection{Proof of Proposition \ref{2}}\label{section_Proof_P2}
To prove the proposition it suffices to integrate the inequality
obtained in Proposition \ref{lemma1} and to use the following lemma
which will be extensively exploited in the sequel.
\begin{lemma}\label{matrixDP} There exists $\lambda>0$ such that
 $\forall n>1$ and $\forall h\in\cH_n$, we have
$$
\lambda_n(h)\geq \lambda.
$$
where $\lambda_n(h)$ is the smallest eigenvalue of the matrix
$$
\cM_{nh}(y)=\frac{1}{nh^{d}}\sum_{i=1}^{n}K^\top\left(\frac{X_i-y}{h}\right)
K\left(\frac{X_i-y}{h}\right)\bI_{V_{h}(y)}(X_i),
$$
and $ K(z) $ is the $D_b$-dimensional vector of the monomials $ z^p,\:p\in\cP_b $.
\end{lemma}

\paragraph{Proof of Proposition \ref{2}.}
In order to simplify the proof, let us introduce the following constants
$$
 c_1=\frac{(1+6D_b^2)}{6A(f)D_b^2},\quad c_2=\frac{\lambda}{432M(f)\:D_b^3}.
$$
By definition of $A(f)$, $M(f)$, $ \mB(.,.) $ respectively in (\ref{borneInfSup}), (\ref{bornepui}) and $ A $, $ M $, we have the following
inequality $ \mB\big(A(f),M(f)\big)\leq \mB(A,M) $. By integration of Proposition \ref{lemma1} and using Lemma \ref{matrixDP}, we get for any $ q\geq1
$ and $ f\in\bH_d(\beta,L,M,A)
$
\begin{eqnarray}
&&\bE_f\big|\hat{f}^h(y)-f(y)\big|^q\bI_G\nonumber\\
&&\quad=\int_0^{+\infty}\eta^{q-1}\bP_f\left(\big|\hat{f}^h(y)-f(y)\big|\geq\eta,G\right)\:d\eta \nonumber\\
&&\quad=(nh^d)^{-q}\int_0^{+\infty}\eta^{q-1}\bP_f\left(\big|\hat{f}^h(y)-f(y)\big|\geq\frac{\eta}{nh^d},G\right)\:d\eta \nonumber\\
&&\quad=(nh^d)^{-q}\left[\int_0^{\frac{2c_1}{c_2}(1\vee\cN_h)}\eta^{q-1}\:d\eta\right.\nonumber\\
&&\qquad\qquad\qquad\left.+\int_{\frac{2c_1}{c_2}(1\vee\cN_h)}^{+\infty}\eta^{q-1}\bP_f\left(\big|\hat{f
} ^h(y)-f(y)\big|\geq\frac{\eta}{nh^d}, G\right)\:d\eta\right] \nonumber\\
&&\quad\leq\frac{(1\vee\cN_h)^q}{(nh^d)^q}\left[\frac{2^qc_1^q}{q\:c_2^q}+\frac{2^q}{c_2^q}\mB\big(A(f),M(f)\big)\Gamma(q)\right]\nonumber,
\end{eqnarray}
where $ \Gamma(\cdot) $ is the well-known Gamma function. 
By definition of $ b_h $ and $ \cN_h $ respectively defined in (\ref{def_bias}) and (\ref{borneexp}), the assertion of Proposition \ref{2}
is proved:
\begin{eqnarray*}
\bE_f\big|\hat{f}^h(y)-f(y)\big|^q\bI_G \leq C^*_q\big(A(f),M(f)\big)\left[\frac{1\vee Ld\:nh^{\beta+d}}{nh^d}\right]^{q},
\end{eqnarray*}
where
$$
C^*_q(a,m)=\frac{1}{q}\left[\frac{432m D_b^3(1+6D_b^2)}{3\lambda aD_b^2}\right]^q+\left[864m\lambda^{-1} D_b^3\right]^q\mB(a,m)\Gamma(q),\quad a,m>0.
$$
\epr
\subsection{Proof of Theorem \ref{minimax}}By definition of $\bar h=(Ln)^{-\frac{1}{\beta+d}}$ and we have
$$
Ldn\bar h^{\beta+d}=d,\quad(n\bar h^d)^{-q}=L^{\frac{qd}{\beta+d}}\varphi_n^q(\beta).
$$ 
Applying the inequality given in Remark \ref{remark_set_norandom}, we come to the assertion of the theorem.
\epr
\subsection{Proof of Theorem \ref{lepski}}This Proof is based on the Lepski scheme developed by \cite{Lepski91} and adapted for the bandwidth
selection by \cite{Lepski_Mammen_Spokoiny97}.
We start the proof with formulating auxiliary Lemmas whose proofs are given in
Appendix (Section \ref{appendixA}). Define 
$$
h^*=\left[n^{-1}c\left(1+\frac{(b-\beta)}{(b+d)(\beta+d)}\ln{n}\right)\right]^{\frac{1}{\beta+d}},
$$
where the positive constant $ c $ is chosen as follows
$$
c<\left[1\wedge 1/(Ld)\right]\left[1\wedge
4/M(f)\right]\frac{\beta+d-1}{\beta+d}\left[1\wedge\frac{A}{144MD_b} \right]\left[1\wedge\frac{6AD_b^2}{1+6D_b^2} \right],
$$
and let the integer $\kappa$ be defined as follows.
\begin{equation}\label{kappa}
2^{-\kappa}h_{\max}\leq h^*<2^{-\kappa+1}h_{\max}.
\end{equation}
The definitions of $h^*$ and $\kappa$ imply the following Lemmas.
\begin{lemma}\label{lem_Proof_Th_lepski}
 \begin{equation*}
\bE_f\big|\hat{f}^{(k)}(y)-f(y)\big|^q\bI_G\leq \bar C_q\frac{(1+k\ln 2)^q}{(nh_k^d)^q},\quad\forall k\geq\kappa,
\end{equation*}
where $ \displaystyle\bar C_q= C_q^*\big(A(f),M(f)\big)\frac{c(\beta+d)}{(\beta+d-1)(Ld)^{-1}}. $
\end{lemma}

\begin{lemma}\label{controlproba}
 For any $f\in
\bH_d(\beta,L,M,A)$ and any   $k\geq \kappa+1$
\begin{eqnarray*}
&&\bP_f\big(\hat{k}=k,\: G\big)\leq J_2\mathfrak{B}(A,M)\exp{\big\{J_1
n(h^*)^{\beta+d}\big\}}2^{-(k-1)(8qd+4)},
\end{eqnarray*}
where $ J_1=Ld(1+6D_b^2)/6AD_b^2 $ and $J_2=(1-2^{-(8qd+4)})^{-1}$.
\end{lemma}
\begin{lemma}\label{eventsAM}
There exists a universal constant $\vartheta>0$ such that
$$
\limsup_{n\rightarrow\infty}\sup_{f\in\bH_d(\beta,L,M,A)}\exp\left\{\frac{An^{\frac{b}{b+d}}}{16M\vartheta^2D_b^2}\right\}\bP_f\big(G^c\big)=0.
$$
\end{lemma}

\paragraph{Proof of Theorem \ref{lepski}.} 
We decompose the risk as follows
\begin{eqnarray}\label{decomposition}
&&\mathbb{E}_f\big|\hat f^{(\hat k)}(y)-f(y)\big|^q\bI_{G}\nonumber\\
&&\quad\leq\mathbb{E}_f\big|\hat f^{(\hat
k)}(y)-f(y)\big|^q\bI_{\hat k\leq \kappa,G}+\mathbb{E}_f\big|\hat
f^{(\hat
k)}(y)-f(y)\big|^q\bI_{\hat k>\kappa,G}\nonumber\\
&&\quad=R_1(f)+R_2(f).
\end{eqnarray}
First we control $R_1$. Obviously
$$
\big|\hat f^{(\hat k)}(y)-f(y)\big|\leq\big|\hat f^{(\hat
k)}(y)-\hat f^{(\kappa)}(y)\big|+\big|\hat
f^{(\kappa)}(y)-f(y)\big|.
$$
Note that  the realization of the event $G$ implies $\hat
M\leq3M(f)/2$. This together with the definition of $\hat k$ yields
$$
\big|\hat f^{(\hat k)}(y)-\hat f^{(\kappa)}(y)\big|\bI_{\hat k\leq
\kappa,G}\leq Cs_n(\kappa),\quad s_n(k)=(1+k\ln2)^q(nh_k^d)^{-q},
$$
where $C=288MD_b^3\lambda^{-1}(32qd+16)$. In view of
Lemma \ref{lem_Proof_Th_lepski} we also get
$$
\bE_f\big|\hat f^{(\kappa)}(y)-f(y)\big|^q\leq\bar C_{q}
s_n(\kappa).
$$
Noting that the right hand side of the obtain inequality is
independent of $f$ and taking into account the definition of
$\kappa$ and $h^*$ we obtain
\begin{eqnarray}\label{r200}
&&\limsup_{n\to\infty}\sup_{f\in\bH_d(\beta,L,A,M)}\phi^{-q}_n(\beta)R_1(f)<\infty.
\end{eqnarray}
Now let us bounded from above $R_2$. Applying Cauchy-Schwartz
inequality we have in view of Lemma \ref{controlproba}
\begin{eqnarray}\label{r2}
R_2(f)&=&\sum_{k>\kappa}^{\mathrm{k}_n}\mathbb{E}_f\big|\hat
f^{(k)}(y)-f(y)\big|^q\bI_{\big[\hat
k=k,G\big]}\nonumber\\
&\leq&\sum_{k>\kappa}\big(\mathbb{E}_f\big|\hat
f^{(k)}(y)-f(y)\big|^{2q}\big)^{1/2}\sqrt{\mathbb P_f\big\{\hat k=k,G\big\}}\nonumber\\
&=&\Delta(h^*)\sum_{k>\kappa}\big(\mathbb{E}_f\big|\hat
f^{(k)}(y)-f(y)\big|^{2q}\big)^{1/2} 2^{-(k-1)(4qd+2)},
\end{eqnarray}
where we have put $\Delta(h^*)=J_2\mathfrak{B}(A,M)\exp{\big\{J_1 n(h^*)^{\beta+d}\big\}}$. We obtain from Lemma \ref{lem_Proof_Th_lepski} and
(\ref{r2})
\begin{eqnarray}
\label{r300}
R_2(f)
\leq \:J_3 \:(nh_{\max}^d)^{-q}\exp{\big\{J_1 n(h^*)^{\beta+d}\big\}},
\end{eqnarray}
where
$$
J_3=J_2\mathfrak{B}(A,M)\:2^{4qd+2}\:\bar{C}_{2q}^{1/2}\sum_{s\geq
0}(1+s\ln 2)^q 2^{-3sdq-2}.
$$
 It remains to note that the definition of $h^*$ implies that
 $$
 \limsup_{n\to\infty}\phi^{-q}_n(\beta)(nh_{\max}^d)^{-q}\exp{\big\{J_1 n(h^*)^{\beta+d}\big\}}<\infty
 $$
and that the right hand side of (\ref{r300}) is  independent of $f$.
Thus,we have
\begin{eqnarray*}
&&\limsup_{n\to\infty}\sup_{f\in\bH_d(\beta,L,A,M)}\phi^{-q}_n(\beta)R_2(f)<\infty.
\end{eqnarray*}
that yields together with (\ref{decomposition}) and (\ref{r200})
\begin{eqnarray*}
&&\limsup_{n\to\infty}\sup_{f\in\bH_d(\beta,L,A,M)}\phi^{-q}_n(\beta)\mathbb{E}_f\big|\hat
f^{(\hat k)}(y)-f(y)\big|^q\bI_{G}<\infty.
\end{eqnarray*}
To get the assertion of the theorem it suffices to show that
\begin{eqnarray}
\label{r400}
\limsup_{n\to\infty}\sup_{f\in\bH_d(\beta,L,A,M)}\phi^{-q}_n(\beta)\mathbb{E}_f\big|\hat
f^{(\hat k)}(y)-f(y)\big|^q\bI_{G^{c}}<\infty.
\end{eqnarray}
Note that  $\hat f^{(\hat k)}(y)\leq 4\hat M$ in view of (\ref{eq:random-parameter-set}). Note also that the local least square estimator
$\tilde\delta$ is linear function of observation $Y^{(n)}$ and, moreover $0\leq Y_i\leq M, i=1,...,n$. This together with the definition of $\hat M$,
(expression (\ref{eq:def-estimators-for-A-M})) allows us to state that there exist $ 0<J_4<+\infty$ such that $\big|\hat f^{(\hat k)}(y)-f(y)\big|\leq
J_4M$. Here we also have taken into account that $||f||_\infty\leq M$.

Finally we obtain
$$
\mathbb{E}_f\big|\hat f^{(\hat k)}(y)-f(y)\big|^q\bI_{G^{c}}\leq
J_4^q\:M^q\bP_f\big\{G^c\big\}.
$$
and (\ref{r400}) follows now from  Lemma \ref{controlproba}.
 \epr
\section{Proofs of lower bounds}\label{sectionlower_bounds}The proofs of Theorems \ref{inf} and \ref{th:optimality-lower-bound}
are based on the following proposition.

Put  $\phi_n(\gamma)=\big[n^{-1}\big(1+(b-\gamma)\ln
n\big)\big]^{\frac{{\gamma}}{{\gamma}+d}},\:\gamma\in (0,b]$ and let
\begin{align*}
R_n^{(q)}(\tilde{f},v)&=\sup_{f\in
\bH_d(\alpha,L,M,A)}\mathbb{E}_f\left[\phi_n^{-q}(\alpha)|\tilde{f}(y)-f(y)|^q\right]\\
&\quad+\sup_{f\in
\bH_d(\beta,L,M,A)}\mathbb{E}_f\left[n^{-vq}\phi_n^{-q}(\beta)|\tilde{f}(y)-f(y)|^q\right].
\end{align*}
where $ v\geq0 $ and $ \alpha,\beta\in(0,b]^2 $.

\begin{proposition}\label{risk}
Let   $\Psi$ be admissible family of normalizations  such that
$$
\psi_n(\alpha)\big/\phi_n(\alpha)\xrightarrow[n\rightarrow\infty]{}0.
$$
Then, for any $0\leq v<(\beta-\alpha)/(\beta+1)(\alpha+1)$
$$
\liminf_{n\rightarrow\infty}\:\inf_{\tilde{f}}\:R_n^{(q)}(\tilde{f},v)>0.
$$
\end{proposition}
The proof is given in section \ref{sectionProof_proposition3}.

\subsection{Proof of Theorem \ref{inf}}
Using the proposition \ref{risk} for $\beta=\alpha$, we have to choose $ v=0 $ and one gets
\begin{eqnarray*}
R_{n,q}\big[\bH_d(\beta,L,M,A)\big]&=&R_n^{(q)}(\tilde{f},0)\\
&=&\sup_{f\in\bH_d(\alpha,L,M,A)}\mathbb{E}_f\left[n^{-q\frac{{\alpha}}{{\alpha}+1}}\big|\tilde{f}(y)-f(y)\big|^q\right]>0,\quad\forall\tilde f.
\end{eqnarray*}\epr
\subsection{Proof of Theorem \ref{th:optimality-lower-bound}}

\textbf{I.} To proof of the first assertion of the theorem it
suffices to consider the family
$\left\{\upsilon_n(\beta)\right\}_{\beta\in (0,b]}$, where
$\upsilon_n(\alpha)=\varphi_n(\alpha)$ and $\upsilon_(\beta)=1$ for
any $\beta\neq\alpha$. The corresponding attainable
estimator is the estimator being minimax on $\bH_d(\alpha,L,M,A).$

\textbf{II.} Let us consider the family
$\left\{\phi_n(\beta)\right\}_{\beta\in ]0,b]}$, which is admissible
in view of Theorem \ref{lepski}. First, we note that $\gamma=b$ is
not possible since $\phi_n(b)=\varphi_n(b)$ the minimax rate of
convergence on $\bH_d(b,L,M,A).$

Thus we  assume that $\gamma$ satisfying
(\ref{eq:th-optimality-lower-bound}) belongs to $\delta\in]0,b[$.
Let $\hat f^{\Psi}$ be a $\Psi^{(n)}$-attainable estimator. Since
$\psi_n(\alpha)/\phi_n(\alpha)\to 0, n\to\infty$ in view of
(\ref{eq:th-optimality-lower-bound}) then obviously
$$
\limsup_{n\rightarrow\infty}\sup_{f\in
\bH_d(\gamma,L,M,A)}\mathbb{E}_f\left[\phi_n^{-q}(\gamma)\big|\hat
f^{\Psi}(y)-f(y)\big|^q\right]=0.
$$
Therefore, applying   Proposition \ref{risk} with $v=0$ we have for
any $\beta<\gamma$
$$
\limsup_{n\rightarrow\infty}\sup_{f\in
\bH_d(\beta,L,M,A)}\mathbb{E}_f\left[\phi_n^{-q}(\beta)\big|\hat
f^{\Psi}(y)-f(y)\big|^q\right]>0.
$$
We conclude that necessarily  $\psi_n(\beta)\gtrsim \phi_n(\beta)$
for any $\beta<\gamma$.

Moreover for any$\beta>\gamma$  applying  Proposition \ref{risk}
with an arbitrary $0\leq v<(\beta-\gamma)/(\beta+1)(\gamma+1)$ we
obtain that
$$
\psi_n(\beta)\gtrsim n^{v}\phi_n(\beta),\:\:\beta>\gamma.
$$
It remains to note that the form of rate of convergence proved in
Theorem \ref{inf} implies that
$$
\phi_n(\gamma)\big/\psi_n(\gamma)=o\left([\ln{n}]^{\frac{\gamma}{\gamma+d}}\right).
$$\epr
\subsection{Proof of Proposition \ref{risk}}\label{sectionProof_proposition3}
Let $\varkappa>0$ the parameter whose choice will be done later. Put
$$
h=\left(\varkappa\frac{1+({\beta}-{\alpha})\ln
n}{n}\right)^{\frac{1}{{\alpha}+d}}.
$$
Later on without loss of generality we will assume  that $L>1$.

Consider the functions: $f_0\equiv 1$ and
\begin{eqnarray*}
f_1(x)=1-\big(L-1\big)\varkappa^{\frac{{\alpha}}{{\alpha}+d}}\phi_n(\alpha)F\left(\frac{x_1-y_1}{h},...,\frac{x_d-y_d}{h}\right),\:
\:\:x\in[0,1]^d.
\end{eqnarray*}
Here F is a compactly  supported positive function belonging to
$\bH_d(\alpha,1,M,A)$  such that $F(0)=1=\max_x F(x)$.

It is easily seen that $f_1\in \bH_d(\alpha,L,M,A)$. Therefore, we
have
\begin{eqnarray*}
R_n^{(q)}(\tilde{f},v)&\geq&\mathbb{E}_0\left|n^{-v}\phi_n^{-1}(\beta)\big(\tilde{f}(y)-1\big)\right|^q+\mathbb{E}_1\left|\phi_n^{-1}
(\alpha)\big(\tilde
{f}(y)-f_1(y)\big)\right|^q\\
&\geq&\mathbb{E}_0\left|n^{-v}\phi_n^{-1}(\beta)\big(\tilde{f}(y)-1\big)\right|^q+\mathbb{E}_1\left|\phi_n^{-1}(\alpha)\big(\tilde{f}
(y)-1\big)+z\right|^q,\\
\end{eqnarray*}
where $z=\big(L-1\big)\varkappa^{\frac{{\alpha}}{{\alpha}+1}}F(0)$.
Set
$$
\tilde{\lambda}=\phi_n^{-1}(\alpha)\big(1-\tilde{f}(y)\big),\quad\varsigma_n=n^{-v}\frac{\phi_n(\alpha)}{\phi_n(\beta)}=n^{-v}\left(\frac{\ln
n}{n}\right)^{-\varrho},
$$
where $\varrho=\frac{{\beta}-{\alpha}}{({\beta}+1)({\alpha}+1)}$. We
get
\begin{eqnarray}
R_n^{(q)}(\tilde{f},v)&\geq&\mathbb{E}_0\big|\varsigma_n\tilde{\lambda}\big|^q+\mathbb{E}_1\big|z-\tilde{\lambda}\big|^q\nonumber\\
&\geq&\mathbb{E}_0\big|\varsigma_n\tilde{\lambda}\big|^q\bI_{\{|\tilde{\lambda}|>z/2\}}+\mathbb{E}_1\big|z-\tilde{\lambda}\big|^q\bI_{\{|\tilde{
\lambda}
|\leq
z/2\}}\nonumber\\
&\geq&\mathbb{E}_0\big|\varsigma_n\frac{z}{2}\big|^q\bI_{\{|\tilde{\lambda}|>z/2\}}+\mathbb{E}_1\big|\frac{z}{2}\big|^q\bI_{\{|\tilde{\lambda}|\leq
z/2\}}\nonumber.
\end{eqnarray}
Noting that  $f_1\leq f_0$, since $F$ is positive, and putting
$c_n\big(Y^{(n)}\big)=\bI_{\{|\tilde{\lambda}|>z/2\}}$ we obtain
\begin{eqnarray}
&&R_n^{(q)}(\tilde{f},v)\geq\varsigma_n^q\frac{z^q}{2^q}\frac{\prod_{i=1}^{n}f_1(X_i)}
{\prod_{i=1}^{n}f_1(X_i)}\int_0^{f_1(X_1)}\ldots\int_0^{f_1(X_{n})}c_n(x)dx_1\ldots
dx_n\nonumber\\
&&\qquad\quad+\frac{z^q}{2^q}\frac{1}{\prod_{i=1}^{n}f_1(X_i)}\int_0^{f_1(X_1)}\ldots\int_0^{f_1(X_{n})}1-c_n(x)dx_1\ldots
dx_n. \label{risque1}
\end{eqnarray}
We have
\begin{eqnarray}
\prod_{i=1}^{n}f_1(X_i)&=&\prod_{i=1}^{n}\left(1-(L-1)\varkappa^{\frac{{\alpha}}{{\alpha}+d}}\phi_n(\alpha)F\left(\frac{X_i-y}{h}
\right)\right)\nonumber\\
&\geq&\left(1-(L-1)\varkappa^{\frac{{\alpha}}{{\alpha}+d}}\phi_n(\alpha)\right)^{nh^d}
\geq e^{-(L-1)\varkappa}n^{-(L-1)\varkappa(\beta-\alpha)}\label{bo}.
\end{eqnarray}
We obtain in view of  (\ref{risque1}) and (\ref{bo})
\begin{eqnarray*}
R_n^{(q)}(\tilde{f},v)&\geq&\varsigma_n^q\frac{z^q}{2^q}e^{-(L-1)\varkappa}n^{-(L-1)\varkappa(\beta-\alpha)}\\
&&\quad\times\frac{1}{\prod_{i=1}^nf_1(X_i)}\int_0^{f_1(X_1)}\ldots\int_0^{f_1(X_{n})}c_n(x)dx_1\ldots
dx_n\\
&&\quad+\frac{z^q}{2^q}\frac{1}{\prod_{i=1}^nf_1(X_i)}\int_0^{f_1(X_1)}\ldots\int_0^{f_1(X_{n})}1-c_n(x)dx_1\ldots
dx_n\\
&\geq&\frac{z^q}{2^q}\left(1\wedge\varsigma_n^qe^{-(L-1)\varkappa}n^{-(L-1)\varkappa(\beta-\alpha)}\right).
\end{eqnarray*}
{\it\underline{Case 1:}} $\beta=\alpha$.  Choosing
$\varkappa=1$, and noting that  $\varsigma_n=1$ and
$\prod_{i=1}^nf_1(X_i)\geq e^{-(L-1)}$, we deduce from
(\ref{risque1}) that yields:
$$
\inf_{\tilde f}R_n^{(q)}(\tilde{f},v)\geq
\frac{(L-1)^q}{2^q}e^{-(L-1)}>0.
$$
{\it\underline{Case 2:}} $\beta>\alpha$.  Put
$$
\varkappa=\frac{q\big(\varrho-v\big)-t_n}{1+(L-1)(\beta-\alpha)}>0,\:\:t_n=\frac{q}{\ln
n}\ln\frac{1}{\big(1+(\beta-\alpha)\ln
n\big)^{-\varrho}}\xrightarrow[n\rightarrow\infty]{}0.
$$
This choice provides us with the following bound
\begin{eqnarray*}
\varsigma_n^qe^{-(L-1)\varkappa}n^{-(L-1)\varkappa(\beta-\alpha)}&=&\big(1+(\beta-\alpha)\ln
n\big)^{-q\varrho}e^{-(L-1)\varkappa}n^{q(\varrho-v)-(L-1)\varkappa(\beta-\alpha)}\\
&\geq&\big(1+(\beta-\alpha)\ln
n\big)^{-q\varrho}e^{-\frac{2}{3}q(L-1)}n^{t_n} \geq
e^{-\frac{2}{3}q(L-1)}.
\end{eqnarray*}
This yields
$$
\inf_{\tilde
f}R_n^{(q)}(\tilde{f},v)\geq\frac{(L-1)^q\varkappa^{q\frac{\alpha}{\alpha+1}}}{2^q}e^{-\frac{2}{3}q(L-1)}>0.
$$
\epr
\section{Appendix}\label{appendixA}

\subsection{Proof of Lemma \ref{propertyZ}}Later on without loss
generality we will suppose that $nh^d\in\bN^*$. In order to simplify understanding of this proof, we note the approximation polynomial $ \cA^i_u=
f_{\theta+u(nh^d)^{-1}}(X_i),\:
i=1,\dots,n $ for all $ u\in \Upsilon_n $.

\medskip

\textbf{1.} Note that for $u\in \Upsilon_n$
\begin{equation}\label{eq_lemme5_EZ}
\bE_fZ_{h,\theta}(u)\leq\prod_{i:\:X_i\in
V_h(y)}\frac{\cA^i_0}{f(X_i)}\leq e^{\cN_h/A(f)}.
\end{equation}
The first inequality is the consequence of the definition of $ Z_{h,\theta} $ in (\ref{not_processZ}) and the following calculation 
$$
\bE_f\bI_{\left[Y_i\leq\:\cA^i_u\right]}=\bP_f\left(Y_i\leq\:\cA^i_u\right)=1\wedge\frac{\cA^i_u}{f (X_i)}.
$$
In (\ref{eq_lemme5_EZ}), the second inequality is obtained with classical inequality $ 1+\rho\leq e^{\rho},\rho\in\bR $ and recall that $
f_\theta(x)\geq f(x) $.
$$
\prod_{i:\:X_i\in V_h(y)}\frac{\cA^i_0}{f(X_i)}=\prod_{i:\:X_i\in
V_h(y)}\left(1+\frac{\cA^i_0-f(X_i)}{f(X_i)}\right)\leq\exp\left\{\frac{b_h\times nh^d}{A(f)}\right\}
$$
\medskip
\\
{\it\underline{Case 1:}} If $\|u_1-u_2\|_1\geq1$, the inequality (\ref{eq_lemme5_EZ}) allows to get
$$
\bE_f\big|Z_{h,\theta}(u_1)-Z_{h,\theta}(u_2)\big|\leq\bE_fZ_{h,\theta}(u_1)+\bE_fZ_{h,\theta}(u_2)\leq2e^{\cN_h/A(f)}\|u_1-u_2\|_1.
$$
\medskip
\\
{\it\underline{Case 2:}} Assume now that $\|u_1-u_2\|_1<1$ and introduce the random events
\begin{eqnarray*}
F_1&=&\left\{\forall i=1,\ldots,n:\:Y_i\leq
\cA^i_{u_1}\wedge \cA^i_{u_2}\right\},\\
F_2&=&\left\{\forall i=1,\ldots,n:\:Y_i\leq
\cA^i_{u_1}\vee
\cA^i_{u_2}\right\}\\
&&\quad\cap\left\{\exists
i:\:Y_i>\cA^i_{u_1}\wedge \cA^i_{u_2}\right\},\\
F_3&=&\left\{\exists i:\:Y_i> \cA^i_{u_1}\vee
\cA^i_{u_2}\right\}.
\end{eqnarray*}
We have used the following notations: $a\wedge b=\min(a,b)$ and
$a\vee b=\max(a,b),\:a,b\in\bR$. For any $(u_1,u_2)\in \Upsilon_n^2$, we
have
\begin{eqnarray}\label{eq_proof_lem1_decomposition}
&&\bE_f\big|Z_{h,\theta}(u_1)-Z_{h,\theta}(u_2)\big|
=\bE_f\big|Z_{h,\theta}(u_1)-Z_{h,\theta}(u_2)\big|\:\bI_{[F_1]}\nonumber\\
&&+\bE_f\big|Z_{h,\theta}(u_1)-Z_{h,\theta}(u_2)\big|\:\bI_{[F_2]}
+\bE_f\big|Z_{h,\theta}(u_1)-Z_{h,\theta}(u_2)\big|\:\bI_{[F_3]}\nonumber\\
&&=\cK_1+\cK_2+\cK_3.
\end{eqnarray}
The following bound will be extensively exploited in the sequel.
$$
f_v(x)\geq2v_{0,...,0}-||v||_1\geq0.25A(f),\quad \forall
v\in\Theta(A(f)/4,9M(f)),\:x\in[0,1]^d.
$$
\paragraph{Control of $\cK_1$.}
\begin{equation}\label{eq_proof_lem1_K1.1}
 \cK_1=\left|\prod_{i:\:X_i\in
V_h(y)}{\frac{\cA^i_0}{\cA^i_{u_1}}}-\prod_{i:\:X_i\in
V_h(y)}{\frac{\cA^i_0}{\cA^i_{u_2}}}\right|\bP_f\big\{F_1\big\},
\end{equation}
and
\begin{equation}\label{eq_proof_lem1_K1.2}
\bP_f\big\{F_1\big\}=\prod_{i:\:X_i\in V_h(y)}\bP_f\big\{Y_i\leq \cA^i_{u_1}\wedge \cA^i_{u_2}\big\}\leq\prod_{i:\:X_i\in
V_h(y)}\frac{ \cA^i_{u_1}\wedge \cA^i_{u_2}}{f(X_i)}.
\end{equation}
Therefore, using (\ref{eq_lemme5_EZ}), we have 
\begin{eqnarray}\label{eq_proof_lem1_K1.3}
\cK_1&\leq&\left(1-\prod_{i:\:X_i\in V_h(y)}\frac{\cA^i_{u_1}\wedge \cA^i_{u_2}}{\cA^i_{u_1}\vee \cA^i_{u_2}}\right)\prod_{i:\:X_i\in
V_h(y)}\frac{\cA^i_0}{f(X_i)}\nonumber\\*[2mm]
&\leq&e^{\cN_h/A(f)}\left(1-\exp\left\{\sum_{i:\:X_i\in V_h(y)}\ln\frac{\cA^i_{u_1}\wedge \cA^i_{u_2}}{\cA^i_{u_1}\vee
\cA^i_{u_2}}\right\}\right).
\end{eqnarray}
Remember that $ \big|\cA^i_{u_1}-\cA^i_{u_2}\big|\leq (nh^d)^{-1}\|u_1-u_2\|_1 $ and $ \cA^i_{u}\geq A(f)/4 $. Let us give the following
calculation with inequality of finite increments for $ \ln(\cdot) $
$$
\ln\frac{\cA^i_{u_1}\wedge \cA^i_{u_2}}{\cA^i_{u_1}\vee \cA^i_{u_2}}
=-\left|\ln \cA^i_{u_1}\wedge \cA^i_{u_2}-\ln \cA^i_{u_1}\vee \cA^i_{u_2}\right|\geq -\frac{(nh^d)^{-1}\|u_1-u_2\|_1}{\cA^i_{u_1}\wedge \cA^i_{u_2}}
$$
Using last inequalities, (\ref{eq_proof_lem1_K1.1}), (\ref{eq_proof_lem1_K1.2}), (\ref{eq_proof_lem1_K1.3}), last inequality and the well known
inequality $
1-e^{-\rho}\leq\rho $, we have
$$
\cK_1\leq \frac{1}{A(f)} e^{\cN_h/A(f)}\|u_1-u_2\|_1.
$$
\paragraph{Control of $\cK_2$.}
We could rewritten
\begin{eqnarray*}
F_2&=&\left\{\forall i=1,\ldots,n:\:Y_i\leq
\cA^i_{u_1}\vee
\cA^i_{u_2}\right\}\\
&&\quad\backslash\left\{\forall i=1,\ldots,n:\:Y_i\leq
\cA^i_{u_1}\wedge
\cA^i_{u_2}\right\}\\
&=&G\backslash F_1.
\end{eqnarray*}
and define
\begin{eqnarray*}
\cG_1&=&\left\{X_i\in V_h(y):\:\cA^i_{u_1}\vee
\cA^i_{u_2}<f(X_i)\right\},\\
\cG_2&=&\left\{X_i\in V_h(y):\:\cA^i_{u_1}\wedge
\cA^i_{u_2}<f(X_i)\right\}.
\end{eqnarray*}
Note that $F_1\subseteq G$ and, therefore,
\begin{eqnarray*}
\cK_2&\leq&\prod_{i:\:X_i\in
V_h(y)}\frac{\cA^i_0}{\cA^i_{u_1}\wedge
\cA^i_{u_2}}\left(\bP_f\big\{G\big\}-\bP_f\big\{F_1\big\}\right)\\
&=&\prod_{i:\:X_i\in
V_h(y)}\frac{\cA^i_0}{\cA^i_{u_1}\wedge
\cA^i_{u_2}}\left(\prod_{i:\:X_i\in \cG_1}\frac{
\cA^i_{u_1}\vee
\cA^i_{u_2}}{f(X_i)}-\prod_{i:\:X_i\in \cG_2}\frac{
\cA^i_{u_1}\wedge \cA^i_{u_2}}{f(X_i)}\right)\\
\end{eqnarray*}
The definition of  $\cG_2$ implies
\begin{eqnarray*}
\prod_{i:\:X_i\in V_h(y)}\frac{1}{\cA^i_{u_1}\wedge \cA^i_{u_2}}\leq\prod_{i:\:X_i\in \cG_2}\frac{1}{\cA^i_{u_1}\wedge \cA^i_{u_2}}\prod_{i:\:X_i\in
\cG_2^c}\frac{1}{f(X_i)}
\end{eqnarray*}
Since $\cG_1\subseteq\cG_2$, $\|u_1-u_2\|_1<1$ and
$|f_u(x)|\leq||u||_1,\,\forall x\in[0,1]^d,\,\forall u\in \Upsilon_n$,
using the last inequality and (\ref{eq_lemme5_EZ}), we obtain
\begin{eqnarray*}
\cK_2&\leq&\prod_{i:\:X_i\in
V_h(y)}\frac{\cA^i_0}{f(X_i)}\left(\prod_{i:\:X_i\in
\cG_2}\frac{\cA^i_{u_1}\vee
\cA^i_{u_2}}{\cA^i_{u_1}\wedge
\cA^i_{u_2}}-1\right)\\
&\leq&4D_b e^{1+\cN_h/A(f)}\|u_1-u_2\|_1/A(f).
\end{eqnarray*}
\paragraph{Control of $\cK_3$.}
We can rewritten the process $ Z_{h,\theta} $ with the notation $ \cA^i_{u} $ 
$$
Z_{h,\theta}(u)=\prod_{i:\:X_i\in
V_h(y)}{\frac{\cA^i_{0}}{\cA^i_{u}}\:\bI_{\left[Y_i\leq\:\cA^i_{u}\right]}}.
$$
Under the event $ F_3 $, we get
$$
\big|Z_{h,\theta}(u_1)-Z_{h,\theta}(u_2)\big|\:\bI_{[F_3]}=0
$$
Then $\cK_3=0$.

\bigskip

The first assertion of the lemma is proved with (\ref{eq_proof_lem1_decomposition}) and the bounds of $ \cK_1 $, $ \cK_2 $ and $ \cK_3 $.

\medskip

\textbf{2.} For any $u\in \Upsilon_n$, since the random variables
$(Y_i)_i$ are independent we have,
$$
\bE_fZ^{1/2}_{h,\theta}(u)=\prod_{i:\:X_i\in
V_h(y)}\sqrt{\frac{\cA^i_0}{\cA^i_u}}\:\bP_f\left\{Y_i\leq
\cA^i_u\right\}.
$$
For any $i$, we have
\begin{eqnarray*}
\sqrt{\frac{\cA^i_0}{\cA^i_u}}\:\bP_f\left\{Y_i\leq\cA^i_u\right\}=\sqrt{\frac{\cA^i_0}{\cA^i_u}}\left[1\wedge\frac{\cA^i_u}{f(X_i)}\right]
\leq\frac{\cA^i_0}{f(X_i)}\left[\frac{f(X_i)}{\sqrt{\cA^i_0}\sqrt{\cA^i_u}}\wedge\frac{\sqrt{\cA^i_u}}{\sqrt{\cA^i_0}}\right].
\end{eqnarray*}
Remind that in view of (\ref{controlbias})
 $f_\theta(x)\geq f(x)$ and $0<f_\theta(x)\leq 3M(f)$ for $x\in
V_h(y)$.

Moreover, for $u\in
\Upsilon_n=nh^d\big(\Theta\big(A(f)/4,9M(f)\big)-\theta\big)$ ,
$0<f_{\theta+u(nh^d)^{-1}}(x)\leq 9M(f)$. Thus for all $i\::\:X_i\in V_h(y)$,
\begin{eqnarray*}
\sqrt{\frac{\cA^i_0}{\cA^i_u}}\:\bP_f\left\{Y_i\leq\cA^i_u\right\}
\leq\frac{\cA^i_0}{f(X_i)}\left[\frac{\sqrt{\cA^i_0}}{\sqrt{\cA^i_u}}\wedge\frac{\sqrt{\cA^i_u}}{\sqrt{\cA^i_0}}\right]
\leq\frac{\cA^i_0}{f(X_i)}\left[1-\frac{\left|\cA^i_0-\cA^i_u\right|}{9M(f)}\right]^{1/2}.
\end{eqnarray*}
The last inequality implies
\begin{eqnarray}\label{final}
\bE_fZ^{1/2}_{h,\theta}(u)&\leq&\prod_{i:\:X_i\in
V_h(y)}\frac{\cA^i_0}{f(X_i)}\sqrt{1-\frac{|f_{u(nh^d)^{-1}}(X_i)|}{9M(f)}}\nonumber\\
&\leq&e^{\cN_h/A(f)}\exp\left\{-\frac{1}{18M(f)\:nh^d}\sum_{i:\:X_i\in
V_h(y)}\big|f_{u}(X_i)\big|\right\}.
\end{eqnarray}
It remains to show
\begin{eqnarray}\label{eq1:new}
\frac{1}{nh^d}\sum_{i:\:X_i\in V_h(y)}\big|f_{u}(X_i)\big|\geq
\lambda_n(h)D_b^{-1}||u||_1.
\end{eqnarray}
Let us remember that $u=(u_{p},\:p\in\cP_b)$ (where $ \cP_b $ is defined in (\ref{polynome})). First, we get from the definition of $f_u$
$$
f_u(x)=u\:K^\top\left(\frac{x-y}{h}\right)=K\left(\frac{x-y}{h}\right)\:u^\top,\quad
\forall x\in[0,1]^d,
$$
and, therefore,
$$
\frac{1}{nh^d}\sum_{i:\:X_i\in
V_h(y)}\big|f_{u}(X_i)\big|=\frac{1}{nh^d}\sum_{i:\:X_i\in
V_h(y)}\left|u\:K^\top\left(\frac{X_i-y}{h}\right)\right|.
$$
Assume $u\neq0$ and put  $v=u/||u||_1$. Noting  that $|f_v(x)|\leq
1,\quad \forall x\in[0,1]^d$, we have
\begin{eqnarray*}
&&\frac{1}{nh^d}\sum_{i:\:X_i\in V_h(y)}\big|f_{u}(X_i)\big|\nonumber\\
&&\quad\geq\frac{1}{nh^d}\sum_{i:\:X_i\in
V_h(y)}\left|u\:K^\top\left(\frac{X_i-y}{h}\right)\right|\:|f_v(X_i)|\nonumber\\
&&\quad=\frac{1}{||u||_1nh^d}\sum_{i:\:X_i\in
V_h(y)}\left|u\:K^\top\left(\frac{X_i-y}{h}\right)K\left(\frac{X_i-y}{h}\right)\:u^\top\right|\nonumber\\
&&\quad\geq\frac{1}{||u||_1}\left|u\:\frac{1}{nh^d}\sum_{i:\:X_i\in
V_h(y)}K^\top\left(\frac{X_i-y}{h}\right)K\left(\frac{X_i-y}{h}\right)\:u^\top\right|.
\end{eqnarray*}
The bound (\ref{eq1:new}) follows now from Lemma \ref{matrixDP}.
The assertion of the lemma follows from (\ref{final}) and (\ref{eq1:new}).

\medskip

\textbf{3.}
In view of Lemma \ref{propertyZ} ({\it 1}), we have
\begin{eqnarray}\label{tempo}
\bE_f\big|Z_{h,\theta}(u)-Z_{h,\theta}(0)\big|\leq{\cC_h\|u\|_1},\quad u\in \Upsilon_n\backslash0.
\end{eqnarray}
Taking into account that $Z_{h,\theta}(0)=1$ we obtain applying
(\ref{tempo}), Fubini's theorem and Tchebychev inequality
\begin{eqnarray*}
&&\mathbb P_f\left\{\int_0^\delta\cdots\int_0^\delta
Z_{h,\theta}(v)dv<\frac{1}{2}\delta^{D_b}\right\}\\
&&\quad=\mathbb P_f\left\{\int_0^\delta\cdots\int_0^\delta
\big(Z_{h,\theta}(v)-Z_{h,\theta}(0)\big)dv<-\frac{1}{2}\delta^{D_b}\right\}\\
&&\quad\leq\mathbb P_f\left\{\int_0^\delta\cdots\int_0^\delta
\big|Z_{h,\theta}(v)-Z_{h,\theta}(0)\big|dv>\frac{1}{2}\delta^{D_b}\right\}\\
&&\quad\leq2\delta^{-D_b}\int_0^\delta\cdots\int_0^\delta\mathbb E_f
\big|Z_{h,\theta}(v)-Z_{h,\theta}(0)\big|dv\\
&&\quad\leq2{\cC_h\delta}
\end{eqnarray*}\epr

\subsection{Proof of Lemma \ref{matrixDP}}{\it\underline{First step:} $ \cM_{{n}h}(y) $ is a nonnegative positive matrix.}\\ Let $\cH_n, n>1$
is defined in
(\ref{Assumptionh}). First, we prove that
\begin{eqnarray}\label{positive_eigen}
\inf_{h\in\cH_n}\lambda_n(h)>0,\quad \forall n>1.
\end{eqnarray}
Suppose that $\exists n_1>1,\:h_{n_1}\in\cH_{n_1}$ such that
$\lambda_{n_1}\big(h_{n_1}\big)=0$. Recall that $ f_t(x)=t\:K(h^{-1}(x-y)) $ for all $ t\in\bR^{D_b} $ and note that $\forall
\tau\in\bR^{D_b}$
\begin{eqnarray*}
\tau^\top\cM_{{n_1}h_{n_1}}(y)\,\tau&=&\frac{1}{nh_{n_1}^d}\sum_{i:\:X_i\in
V_{h_{n_1}}(y)}\left[\tau\:K^\top\left(\frac{X_i-y}{h_{n_1}}\right) \right]^2\\
&=&\frac{1}{nh_{n_1}^d}\sum_{i:\:X_i\in V_{h_{n_1}}(y)}\big[f_\tau(X_i)\big]^2\geq0.
\end{eqnarray*}
Since $\lambda_{n_1}\big(h_{n_1}\big)$ is the smallest eigenvalue of
the matrix $\cM_{{n_1}h_{n_1}}(y)$ the assumption
$\lambda_{n_1}\big(h_{n_1}\big)=0$ implies that there exist $\tau^*$
belonging to the unit sphere of $\bR^{D_b}$ such that
$$
\frac{1}{nh_{n_1}^d}\sum_{i:\:X_i\in
V_{h_{n_1}}(y)}\big[f_\tau(X_i)\big]^2=0.
$$
It obviously implies that $f_{\tau^*}(X_i)=0$ for all $X_i\in
V_{h_{n_1}}(y)$. It remains to note that
$nh_{n_1}^d\geq\big(b+1\big)^d$ since $h_{n_1}\in\cH_n$ and to
apply the result obtained in \cite{Nemirovski00} (page 20). It yields
$\tau^*=0$ and  the obtained contradiction proves
(\ref{positive_eigen}).

\medskip

{\it\underline{Second step:} $ \cM_{nh}(y)\xrightarrow[]{n\rightarrow\infty}\cM $.}\\ Let $\lambda_0$ be the smallest eigenvalue of the
matrix
$$
\cM=\int_{[-1/2,1/2]^d}K^\top (x)\:K(x)\:dx
$$
whose  general term is given by
$$
\cM_{p,q}=\prod_{j=1}^d\int_{-\frac{1}{2}}^{\frac{1}{2}}x_j^{p_j+q_j}dx_j,\quad
0\leq|p|,|q|\leq b.
$$
Let us prove that
\begin{eqnarray}\label{limit_eigen}
\displaystyle
\limsup_{n\rightarrow\infty}\sup_{h\in\cH_n}\big|\lambda_n(h)-\lambda_0\big|=0.
\end{eqnarray}
Put $m=n^{1/d}$ and without loss of generality we will assume that
$m$ is integer. Remind that the general term of the matrix
$\cM_{nh}(y)$ is given by
 $$
\big(\cM_{nh}(y)\big)_{p,q}=\frac{1}{nh^d}\sum_{i:\:X_i\in
V_h(y)}\prod_{j=1}^d\left(\frac{X_{i_j}-y_j}{h}\right)^{p_j+q_j}.
$$
where $X_{i_j}=i_j/m$ for all $ j=1,\dots,d $ and $X_i=\big(X_{i_1},...,X_{i_d}\big)$. We get
\begin{eqnarray*}
&&\frac{1}{nh^d}\sum_{i:\:X_i\in
V_h(y)}\prod_{j=1}^d\int_{i_j-1}^{i_j}\left(\frac{x_j/m-y_j}{h}\right)^{p_j+q_j}dx_j\\
&&\quad\leq\frac{1}{nh^d}\sum_{i:\:X_i\in
V_h(y)}\prod_{j=1}^d\left(\frac{X_{i_j}-y_j}{h}\right)^{p_j+q_j}\\
&&\quad\leq\frac{1}{nh^d}\sum_{i:\:X_i\in
V_h(y)}\prod_{j=1}^d\int_{i_j}^{i_j+1}\left(\frac{x_j/m-y_j}{h}\right)^{p_j+q_j}dx_j,
\end{eqnarray*}
It yields by change of variables that
\begin{eqnarray}
\label{eq10:new}
\prod_{j=1}^d\int_{-\frac{1}{2}-2(nh^d)^{-1}}^{\frac{1}{2}}x_j^{p_j+q_j}dx_j&\leq&\frac{1}{nh^d}\sum_{i:\:X_i\in
V_h(y)}\prod_{j=1}^d\left(\frac{X_{i_j}-y_j}{h}\right)^{p_j+q_j}\nonumber\\
&\leq&\prod_{j=1}^d\int_{-\frac{1}{2}}^{\frac{1}{2}+2(nh^d)^{-1}}x_j^{p_j+q_j}dx_j,
\end{eqnarray}
Note that  $\displaystyle nh^d\geq\ln^{\frac{1}{1+d}}(n)$ for any
$h\in\cH_n$. This together with (\ref{eq10:new}) yields
$$
\limsup_{n\rightarrow\infty}\sup_{h\in\cH_n}\left|\big(\cM_{nh}(y)\big)_{p,q}-\cM_{p,q}\right|=0,\quad 0\leq|p|,|q|\leq b.
$$
The last result obviously imply (\ref{limit_eigen}).

\medskip

{\it\underline{Third step:} Conclusion.}\\ First we show that $\lambda_0>0$. Indeed, $\forall
\tau\in\bR^{D_b}$
$$
\tau^\top\cM\tau=\int_{[-1/2,1/2]^d}\big[f_\tau(x)\big]^2dx\geq0.
$$
Since  $\lambda_0$ is smallest eigenvalue of the matrix $\cM$ the
assumption $\lambda_0=0$ would imply that there exists $\tau^*$
belonging to the unit sphere of $\bR^{D_b}$ such that
$f_{\tau^*}\equiv0$.  Since $f_{\tau^*}$ is a polynomial the last
identity is possible if and only if $\tau^*=0$. The obtained
contradiction shows that $\lambda_0>0$.

Next, note  that in view of (\ref{limit_eigen}) there exists $n_0$
such that $\forall n>n_0$ and $\forall h\in\cH_n,$ $
\lambda_n(h)\geq\lambda_0/2. $

On the other hand in view of (\ref{positive_eigen}) $ \min_{n\leq
n_0}\inf_{h\in\cH_n}\lambda_n(h)>0. $ It remains to define
$\lambda>0$ as
$$\displaystyle\lambda=\min\left(\min_{n\leq
n_0}\inf_{h\in\cH_n}\lambda_n(h),\:\lambda_0/2\right).
$$
\epr
\subsection{Proof of Lemma \ref{lem_Proof_Th_lepski}}
Remind that $ h_k\leq h_\kappa\leq h^* $ by definition of $ h_k $, $ h^* $ and $ \kappa $ (see (\ref{kappa})). Using Proposition \ref{2} with $ h=h_k
$, it yields 
 \begin{eqnarray}\label{eq_proof_lem3.1}
\bE_f\big|\hat{f}^{(k)}(y)-f(y)\big|^q\bI_G&\leq& C_q^*\big(A(f),M(f)\big)\left(\frac{1\vee Ld\:nh_k^{\beta+d}}{nh_k^d}\right)^q\nonumber\\
&\leq& C_q^*\big(A(f),M(f)\big)\left(\frac{1\vee Ld\:n(h^*)^{\beta+d}}{nh_k^d}\right)^q.
\end{eqnarray}
The control of $ n(h^*)^{\beta+d} $ requires the following calculation.
\begin{equation}\label{eq_proof_lem3.2}
n(h^*)^{\beta+d}\leq 1+\frac{b-\beta}{(b+d)(\beta+d)}\ln n=\rho_n(\beta)
\end{equation}
where $ \rho_n(\beta) $ is the price to pay for adaptation defined in (\ref{adaptive_rate}). By definition of $ h_k $, we have
\begin{eqnarray*}
 1+\kappa\ln2&=&1+\ln\frac{h_{\max}}{h_k}\geq1+\ln\frac{h_{\max}}{h*}\\
  &\geq&1+\frac{b-\beta}{(b+d)(\beta+d)}\ln
n-\frac{1}{\beta+d}\ln\big[c\big(1+(b-\beta)\ln n\big)\big].
\end{eqnarray*}
Using the classical inequality $ \ln(1+x)\leq x $ and $ c\leq 1 $, we obtain with the last inequality
\begin{equation}\label{eq_proof_lem3.3}
\frac{\beta+d-1}{\beta+d}\rho_n(\beta)\leq1+\kappa\ln2\leq1+k\ln2, \forall k\geq\kappa.
\end{equation}
According to (\ref{eq_proof_lem3.1}), (\ref{eq_proof_lem3.2}) and (\ref{eq_proof_lem3.3}), Lemma \ref{lem_Proof_Th_lepski} is proved.
\epr

\subsection{Proof of Lemma \ref{controlproba}}
Note that for any $k\geq \kappa+1$ and by definition of $\hat k$ in (\ref{indexe adaptive})
$$
\big\{\hat{k}=k\big\}=\cup_{l\geq k}\left\{\big|\hat
f^{(k-1)}(y)-\hat f^{(l)}(y)\big|> \hat{M} S_n(l)\right\}.
$$
Note that $S_n(l)$ is monotonically increasing in $l$ and,
therefore,
\begin{eqnarray*}
\big\{\hat{k}=k\big\}&\subseteq &\left\{\big|\hat
f^{(k-1)}(y)-f(y)\big|> 2^{-1}\hat{M} S_n(k-1)\right\}
\\*[2mm]
& &\cup\left[ \cup_{l\geq k}\left\{\big|\hat f^{(l)}(y)-f(y)\big|>
2^{-1}\hat{M} S_n(l)\right\}\right].
\end{eqnarray*}
Taking into account that the event $G$ implies the realization of
the event $\hat M\geq M(f)/2\geq A/2$ we come to the following
inequality: for any $k\geq \kappa+1$
\begin{eqnarray}
\label{eq1:proof-lemma2} \bP\big(\hat{k}=k,\:
G\big)&\leq&\bP\left\{\big|\hat f^{(k-1)}(y)-\hat f(y)\big|>
4^{-1}M(f)\: S_n(k-1),\: G\right\} \nonumber\\*[2mm] & &+
\sum_{l\geq k}\bP\left\{\big|\hat f^{(l)}(y)-f(y)\big|> 4^{-1}M(f)\:
S_n(l),\: G\right\}.
\end{eqnarray}
Now we go to justify the use of Proposition \ref{lemma1}. Note that $b_{h_l}\leq Ld h_l^{\beta}$ since
$f\in\bH_d(\beta,L,A,M)$ and, therefore, by definition of $ h^* $, we have
\begin{eqnarray}
\label{eq3:proof-lemma2}\cN_{h_l}\leq Ldn(h_l)^{\beta+d}\leq
Ldn(h_\kappa)^{\beta+d}\leq Ldn(h^*)^{\beta+d}\leq c\rho_n(\beta)  , \:\:\forall l\geq k-1.
\end{eqnarray}
Remark that the definition of $S_n(l)$ yields
$$
nh_l^d\:S_n(l)\geq
432D_b^3(32qd+16)\lambda^{-1}(h_l)\big[1+\ln{\big(h_{\max}/h_l}\big)\big].
$$
Using (\ref{eq_proof_lem3.3}), (\ref{eq3:proof-lemma2}) and the last inequality, we have 
\begin{equation}\label{eq4:proof-lemma2}
\frac{M(f)}{4} nh_l^d\:S_n(l)\geq 144MD_b(1\vee\cN_{h_l})/(\lambda_n(h_l) A).
\end{equation}
 
The last inequality allows us to apply Proposition \ref{lemma1} and Lemma \ref{matrixDP}
with $\ve=\frac{M(f)}{4} nh_l^d\:S_n(l)$,
we obtain $\forall l\geq k-1$
\begin{eqnarray}
\label{eq2:proof-lemma2} &&\bP\left\{\big|\hat f^{(l)}(y)-f(y)\big|>
(M(f)/4)\: S_n(l),\: G\right\}\nonumber\\
&&\quad\leq
\mathfrak{B}(A,M)\cE(h_l)\left[h_{\max}/h_l\right]^{-8qd-4}
\nonumber\\
&&\quad=\mathfrak{B}(A,M)\cE(h_l)2^{-l(8qd+4)}.
\end{eqnarray}

Here we have also used that $k\geq \kappa+1$. We obtain from
(\ref{eq1:proof-lemma2}),  (\ref{eq2:proof-lemma2}) and
(\ref{eq3:proof-lemma2}) that $k\geq \kappa+1$
\begin{eqnarray*}
&&\bP\big(\hat{k}=k,\: G\big)\leq J_2\mathfrak{B}(A,M)\exp{\big\{J_1
n(h^*)^{\beta+d}\big\}}2^{-(k-1)(8qd+4)},
\end{eqnarray*}
where $J_2=(1-2^{-(8qd+4)})^{-1}$. \epr

\subsection{Proof of Lemma \ref{eventsAM}}
Put for any $p\in \cP_b$
$$
W_{ni}^{p}(y)=p_1!...p_d!\:\frac{h_{\max}^{d-|p|}}{n}\:K^\top(0)\:\cM_{nh_{\max}}^{-1}(y)\:
K\left(\frac{X_i-y}{h_{\max}}\right)\bI_{V_{\max}(y)}(X_i),
$$
and note that $ \tilde\delta_{p}=\sum_{i=1}^d 2 Y_i\: W_{ni}^{p}(y). $

The model (\ref{model}) can be rewritten as
$2Y_i=f(X_i)+f(X_i)(2U_i-1)$. Thus, putting
$F(X)=\big(f(X_i)\big)_{i=1,...,n}$,
$V(X)=\big(f(X_i)(2U_i-1)\big)_{i=1,...,n}$ and
\[
\cD(f)=\left(\displaystyle \frac{\partial^{|p|}f(y)}{\partial
y_1^{p_1}\cdots\partial y_d^{p_d}},\:p\in\cP_\beta\right),
\]

 $1^0$. {\it Deviations of $ \hat M $.} By definition of $ \hat M $ in (\ref{eq:def-estimators-for-A-M}), we obtain
\[
|\hat
M-M(f)|\leq||\tilde\delta-\cD(f)||_1\leq\left\|\cV\:F(X)-\cD(f)\right\|_1+\left\|\cV\:V(X)\right\|_1.
\]
Here $\cV$ is $D_b\times n$-matrix  of general term $ \cV_{p
i}=W_{ni}^{p}(y) $ and $ \|.\|_1 $ is the $ \ell_1 $-norm. Let us prove that
\begin{equation}\label{eq:proof4_inter}
\bP_f\left\{|\hat M-M(f)|>M(f)/2\right\}\leq\exp\left\{-\frac{n^{\frac{b}{b+d}}}{8\vartheta_2^2D_b^2}\right\}.
\end{equation}
In view of the result proved in \cite{Hardle_Hart_Marron_Tsybakov92} and \cite{Tsybakov08} there exist $\vartheta_1,\vartheta_2>0$ such
that
\begin{eqnarray*}
\left\|\cV\:F(X)-\cD(f)\right\|_1&\leq&\vartheta_1h_{\max}^{\beta-\lfloor\beta\rfloor},\\
\sup_{i,x}|W_{ni}^{p}(y)|&\leq&\frac{\vartheta_2}{nh_{\max}^d},\quad
p\in\cP_\beta.
\end{eqnarray*}
Remind that $h_{\max}\xrightarrow[]{n\rightarrow\infty}0$ and, therefore, $\exists n_0$ such that $\vartheta_1h_{\max}^{\beta-\lfloor\beta\rfloor}\leq
M(f)/4$ for any $n\geq n_0$. Note
that $n_0$ can be chosen independent on $f$ since $M(f)/4\geq A/4$.
Thus, we get
\begin{eqnarray*}
&&\bP_f\left\{|\hat M-M(f)|>M(f)/2\right\}\\
&&\quad\leq\sum_{p\in\bN^d:\,0\leq|p|\leq\beta}\bP_f\left\{\left|\sum_{X_i\in[0,1]^d}f(X_i)(2U_i-1)W_{ni}^{p}(y)\right|>\frac{M(f)}{4D_b}\right\}.
\end{eqnarray*}
Noting that $ \displaystyle\big|f(X_i)(2U_i-1)W_{ni}^{p}(y)\big|\leq M(f)\frac{\vartheta_2}{nh_{\max}^d}, $ applying Höeffding inequality
\cite{Boucheron_Bousquet_Lugosi04} and the last inequality, we obtain
\begin{eqnarray}
\label{eq2:new} &&\sum_{p\in\bN^d:\,0\leq|p|\leq
\beta}\bP_f\left\{\left|\sum_{X_i\in[0,1]^d}f(X_i)(2U_i-1)W_{ni}^{p}(y)\right|>\frac{M(f)}{4D_b}\right\}\nonumber\\
&&\quad\leq
D_b\exp\left\{-\frac{nh_{\max}^d}{8\vartheta_2^2D_b^2}\right\}=D_b\exp\left\{-\frac{n^{\frac{b}{b+d}}}{8\vartheta_2^2D_b^2}\right\}.
\end{eqnarray}
Therefore (\ref{eq:proof4_inter}) is proved. 

 $2^0$. {\it Deviations of $ \hat A $.} Since $|f(y)-A(f)|\leq Ldh_{\max}^\beta\leq A(f)/4$ for $n\geq n_0$ one
has
$$
\bP_f\left\{|\hat A-A(f)|>A(f)/2\right\}\leq\bP_f\left\{|\hat
M-M(f)|>A(f)/4\right\}.
$$
Repeating previous calculations we obtain
\begin{eqnarray}
\label{eq3:new} \bP_f\left\{|\hat A-A(f)|>A(f)/2\right\}&\leq&
D_b\exp\left\{-\frac{\big[A(f)\big]^2n^{\frac{b}{b+d}}}{16\big[M(f)\big]^2\vartheta_2^2D_b^2}\right\}\nonumber\\
&\leq&
D_b\exp\left\{-\frac{An^{\frac{b}{b+d}}}{16M\vartheta_2^2D_b^2}\right\}.
\end{eqnarray}
Since
$\bP_f\big(G^c\big)\leq\bP_f\big(G_{\hat A}^c\big)+\bP_f\big(G_{\hat M}^c\big)$
the assertion of the lemma follows from (\ref{eq2:new}) and
(\ref{eq3:new}). \epr

\section*{Acknowledgements}
I would like to thank Oleg Lepski for his helpful remarks and comments helped to improve the presentation of the paper. 
I am extremely grateful to two referees and the Associate Editor for very helpful comments which are greatly appreciated.

{
\bibliographystyle{plainnat}
\bibliography{reference2}
}

\end{document}